\documentclass[12pt]{amsart}

\date{June 13, 2015}

\usepackage{amsmath,amsthm,amsfonts,amssymb,multicol}
\usepackage{graphicx}
\usepackage{booktabs}
\usepackage{array}
\usepackage{subfigure}
\usepackage{enumerate}
\usepackage{url}

\usepackage[english]{babel}
\usepackage{verbatim} 
\usepackage{enumitem}
\usepackage{bbm} 
\usepackage[normalem]{ulem} 
\usepackage{bm}
\usepackage{mathtools}
\usepackage{hyperref}
\usepackage{graphicx}

\newcommand{\be}{\begin}
\newcommand{\e}{\end}
\newcommand{\beq}{\begin{equation}}
\newcommand{\eeq}{\end{equation}}
\newcommand{\beqs}{\begin{equation*}}
\newcommand{\eeqs}{\end{equation*}}
\renewcommand{\l}{\left}
\renewcommand{\r}{\right}

\renewcommand{\d}{\mathrm{d}} 
\newcommand{\G}{\curly{G}}

\newcommand{\set}[1]{\mathbb{#1}}

\newcommand{\curly}[1]{\mathcal{#1}}

\newcommand{\setof}[2]{\left\{ #1\; : \;#2 \right\}}

\newcommand{\Om}{\Omega}
\newcommand{\eps}{\epsilon}

\newcommand{\lam}{\lambda}

\newcommand{\blam}{\mathbf{\lambda}}
\newcommand{\al}{\alpha}
\newcommand{\de}{\delta}
\newcommand{\Del}{\Delta}

\newcommand{\ol}{\overline}
\newcommand{\ci}{\circ}

\newcommand{\ind}{\mathbbm{1}}		

\newcommand{\bd}{\mathbf{d}}
\newcommand{\bg}{\mathbf{g}}

\renewcommand{\it}{\infty}

\newcommand{\del}{\partial}

\newcommand{\n}{\mathrm{in}}			
\newcommand{\out}{\mathrm{out}}

\theoremstyle{definition}

\theoremstyle{remark}

\def\dotuline{\bgroup
  \ifdim\ULdepth=\maxdimen  
   \settodepth\ULdepth{(j}\advance\ULdepth.4pt\fi
  \markoverwith{\begingroup
  \advance\ULdepth0.08ex
  \lower\ULdepth\hbox{\kern.15em .\kern.1em}%
  \endgroup}\ULon}

\def\dashuline{\bgroup
  \ifdim\ULdepth=\maxdimen  
   \settodepth\ULdepth{(j}\advance\ULdepth.4pt\fi
  \markoverwith{\kern.15em
  \vtop{\kern\ULdepth \hrule width .3em}%
  \kern.15em}\ULon}
\allowdisplaybreaks

\title{Heat flows on hyperbolic spaces}

\author[M. Lemm]{Marius Lemm}
\address{Mathematics Dept. MC 253-37, Caltech, Pasadena, CA 91125}
\email{mlemm@caltech.edu}

\author[V. Markovic]{Vladimir Markovic}
\address{Mathematics Dept. MC 253-37, Caltech, Pasadena, CA 91125}
\email{markovic@caltech.edu}

\thanks{Vladimir Markovic is supported by the NSF grant number DMS-1500951}

\begin{document}

\begin{abstract}  In this paper we develop new methods for studying the convergence problem for the heat flow on negatively curved spaces and prove that any quasiconformal map of the sphere $\set{S}^{n-1}$, $n\geq 3$, can be extended to the $n$-dimensional hyperbolic space such that the heat flow starting with this extension converges to a quasi-isometric harmonic map. This implies the Schoen-Li-Wang conjecture that every quasiconformal map of $\set{S}^{n-1}$, $n\geq 3$, can be extended to a harmonic quasi-isometry of the $n$-dimensional hyperbolic space. 
\end{abstract}

\maketitle

\setcounter{tocdepth}{1}
\tableofcontents

\section{Introduction and main result}
\subsection{Harmonic maps via heat flows}
A central question in the theory of harmonic maps is under what assumptions a map $\phi:M\rightarrow N$ between two negatively curved Riemannian manifolds can be  deformed into a harmonic map. 

In the pioneering work of Eells and Sampson \cite{EellsSampson64}, it was proved that any $C^1$ map $\phi: M\rightarrow N$ can be deformed into a harmonic map when $M$ and $N$ are compact without boundary and $N$ has negative curvature. Their seminal idea was to obtain the harmonic map as the large time limit of a solution to the heat equation
\beq
\label{eq:HE}
\begin{aligned}
\tau(u)(x,t) &= \del_t u(x,t),\quad \textnormal{on } M\times [0,\it)\\
u(x,0)&=\phi(x),\quad \textnormal{on } M.
\end{aligned}
\eeq
Here $\tau$ denotes the tension field of a map. The convergence of the heat flow as $t\rightarrow\it$ is based on the fact that there is a monotone decreasing energy functional. Importantly, this energy is finite for all initial $C^1$ maps in the compact setting. Hamilton \cite{Hamilton75} proved an analogous statement for compact manifolds with boundary. When $M,N$ are noncompact, convergence of the heat flow was established by Liao and Tam \cite{LiaoTam} under the assumption that $\phi$ has finite energy (see also \cite{SchoenYau76}).  Li and Tam \cite{LiTam91} proved convergence to a harmonic map assuming that $\tau(\phi)\in L^p$ for some $1<p<\it$, see also \cite{Li}. Wang \cite{Wang98} showed that it is enough to assume that $|\tau(\phi)|$  tends to zero uniformly near the boundary to make sure that the heat flow converges. We refer the reader to \cite{Wang98, LinWang} for further background on the heat equation on Riemannian manifolds. We note that the existence of harmonic maps can also be proved without using the heat equation, see e.g.\ \cite{SchoenYau76}.

\subsection{The Schoen-Li-Wang conjecture}
Of particular interest is the case where $M=N=\set{H}^n$ is the $n$-dimensional hyperbolic space. The homotopy class of a map $\phi:\set{H}^n\rightarrow\set{H}^n$ corresponds to its action on the ``boundary'', which we identify as usual with $\set{S}^{n-1}$. The main conjecture is that any quasiconformal boundary map gives rise to a harmonic map of hyperbolic space.

\be{conj}[Schoen, Li and Wang]
\label{conj:schoen}
Let $n\geq 2$. For every quasiconformal map $f:\set{S}^{n-1}\rightarrow \set{S}^{n-1}$, there exists a unique harmonic and quasi-isometric extension $\curly{H}(f):\set{H}^n\rightarrow \set{H}^n$.
\e{conj}

The precise definitions will be given later. Schoen \cite{Schoen93} conjectured this for $n=2$ and the generalization to all $n\geq 2$ is due to Li and Wang \cite{LiWang98}. The uniqueness part of the conjecture was established by Li and Tam \cite{LiTam93II} for $n=2$ and by Li and Wang \cite{LiWang98} for all $n$. The existence part remained an open problem for all $n\geq 2$, with several partial and related results \cite{LiTam93, TamWan98, HardtWolf97, Wang98, BonsanteSchlenker}.  Recently, existence was proved in $n=3$ \cite{Markovic14}  (without using the heat flow). The proof of existence in the $n=2$ case has been announced in \cite{Markovic15}.

Many of the convergence results for the heat equation that we discussed above were motivated by versions of Conjecture \ref{conj:schoen}. The idea is that starting from a quasiconformal boundary map $f:\set{S}^{n-1}\rightarrow \set{S}^{n-1}$, one defines an appropriate extension to hyperbolic space $\set{H}^n$. If the extension is sufficiently regular (e.g.\ has tension field in $L^p$ for some $1<p<\it$), one can run a heat flow with it as the initial map, which then converges to a harmonic map. Since the heat flow will be a quasi-isometry with uniformly bounded distance from the inital map, it is also an extension of $f$. This yields the existence of a harmonic extension if one has a sufficiently regular extension of the quasiconformal boundary map.

The limitations of previous works with regards to the general Conjecture \ref{conj:schoen} lie in the fact that in order to get sufficient regularity of the extension, one needs much stronger regularity of the boundary map $f:\set{S}^{n-1}\rightarrow \set{S}^{n-1}$ than just quasiconformality (it is required that $f$ is $C^1$).  For this reason, the heat flow method has not been successful in proving Conjecture \ref{conj:schoen} so far.

\subsection{Main result}
In this paper, we prove that any quasiconformal map has a ``good extension'' such that the heat flow starting with this extension converges to a harmonic quasi-isometry. Moreover, the regularity of the harmonic map depends only on the distortion $K$ of the quasiconformal map.

\be{thm}
\label{thm:conv}
Let $n\geq 3$ and $K\geq 1$. Let $f:\set{S}^{n-1}\rightarrow \set{S}^{n-1}$ be a $K$-quasiconformal map. Then, there exists a quasi-isometric extension of $f$, $\curly{E}(f):\set{H}^n\rightarrow \set{H}^n$, such that the solution $u(x,t)$ to \eqref{eq:HE} with the choice $\phi\equiv \curly{E}(f)$ converges to a harmonic quasi-isometry $\curly{H}(f)$. 

Moreover, there exist $L=L(K)>0$ and $A=A(K)\geq0$ such that $\curly{E}(f)$ and $\curly{H}(f)$ are $(L,A)$-quasi-isometries.
\e{thm}

\be{rmk}
Throughout the paper we write $C=C(K_1,K_2,\ldots)$ to say that the constant $C$ depends only on $K_1,K_2,\ldots$ The constant $C$ may also implicitly depend on the dimension $n$.
\e{rmk}

The extension $\curly{E}$ is a higher-dimensional generalization of the ``good extension'' constructed recently in \cite{Markovic14}, see Section \ref{sect:good} for the details. We note

\be{cor}
\label{cor:schoen}
Let $n\geq 3$. For every quasiconformal map $f:\set{S}^{n-1}\rightarrow \set{S}^{n-1}$, there exists a harmonic and quasi-isometric extension $\curly{H}(f):\set{H}^n\rightarrow \set{H}^n$.
\e{cor}

Together with the uniqueness result of \cite{LiWang98}, this proves Conjecture \ref{conj:schoen} when $n\geq 3$.

\subsection{A sketch of the proof}

Let $a\in\set{S}^{n-1}$ and write $\curly{G}_a(f)$ for the good extension of a $K$-quasiconformal map $f:\set{S}^{n-1}\rightarrow \set{S}^{n-1}$ for which $f(a)=a$ (see Section \ref{sect:good} for details). Most importantly, the construction is such that $|\tau(\curly{G}_a(f))|$ is small at a ``random'' point in hyperbolic space (i.e.\ the fraction of points on any geodesic sphere where the tension field is greater than $\eps$ goes to zero as the geodesic distance increases, for every $\eps>0$).

We write $u_{a,f}(x,t)$ for the solution to the heat equation \eqref{eq:HE} with initial map $\phi\equiv \curly{G}_a(f)$. (It follows from standard results about the heat equation that $u_{a,f}(x,t)$ exists for all times and is unique, see Proposition \ref{prop:HE}.) The proof of \textbf{Theorem \ref{thm:conv}} is based on the following \emph{two key results}.

For a function $g$ defined on $\set{H}^n$, we write
\beqs
\|g\|=\sup_{x\in \set{H}^n} |g(x)|.
\eeqs

\be{enumerate}
\item[(I)] In \textbf{Theorem \ref{thm:heat}}, we prove $\lim_{t\rightarrow\it}\|\tau(u_{a,f})(\cdot,t)\|=0$.
\item[(II)] In \textbf{Theorem \ref{thm:close}}, we show that there exist $\eps_0=\eps_0(K)>0$ and $D_0=D_0(K,a)$ such that
\beqs
\|\d_{\set{H}^n}\l(\curly{G}_a(f),\psi\r)\|\leq D_0.
\eeqs
holds for all $C^2$ maps $\psi:\set{H}^n\rightarrow\set{H}^n$ which extend $f$ and satisfy $\|\tau(\psi)\|<\eps_0$
\e{enumerate}
Together, Theorems \ref{thm:heat} and \ref{thm:close} readily imply Corollary \ref{cor:schoen}. We will give the complete proof in the next section. Here we just note that by combining them one gets
\beqs
\sup_{t>0}\|\d_{\set{H}^n}\l(u_{a,f}(\cdot,t), \curly{G}_a(f)\r)\|<\it,
\eeqs
which by Arzela-Ascoli and Theorem \ref{thm:heat} implies that $u_{a,f}$ converges to a harmonic map along a subsequence of times $t_i\rightarrow\it$. The limit is still an extension of $f$ because it is a quasi-isometry which is at finite distance from the quasi-isometry $\curly{G}_a(f)$ extending $f$. Uniqueness of the limit then gives convergence for all $t\rightarrow\it$.

While Theorem \ref{thm:close}, will follow by essentially a straightforward generalization of the arguments in \cite{Markovic14}, Theorem \ref{thm:heat} requires more work. It is based on three ingredients: (a) Hamilton's parabolic maximum principle \eqref{eq:tauestimate} for subsolutions of the heat equation, (b) the diffusion of heat in hyperbolic space (see Appendix \ref{app:heat}) and (c) the new Sector Lemma \ref{lm:sector} for the good extension.

\subsection{Discussion}
The theory of good extensions of quasiconformal maps was initiated in \cite{Markovic14}. The most important property of any good extension is that it is ``almost harmonic'' (i.e.\ it has small tension field at ``most'' points). However, in this paper we have to develop  a broader and more detailed theory of  the good extension than the one defined in \cite{Markovic14}. We introduce a family of good extensions $\{\G_a\}_a$ indexed by boundary points at which they are ``anchored'', see Definition \ref{defn:gadefn}. We extend the theory of good extensions by the new notion of \emph{``partial conformal naturality''}. It is important to relate different members of the family $\{\G_a\}_a$. Indeed, it says that for two points $a,b\in \set{S}^{n-1}$ and $I,J\in \mathrm{Isom}(\set{H}^n)$ with $I(b)=J(b)=a$, the good extension ``anchored'' at $a$ and $b$ are related by
\beqs
I\ci \curly{G}_b(f)\circ J^{-1}= \curly{G}_a(I \circ f\circ J^{-1}),
\eeqs
for every $f\in\mathbf{QC}_b(\set{S}^{n-1})$. In particular, this implies that the good extension $\G_a$ is continuous in $a$. We refer to Section \ref{sect:good} for a detailed discussion of the good extension.

We conclude the introduction with the following two remarks.

The recent work \cite{Markovic14}, which proves the existence part of Conjecture \ref{conj:schoen} when $n=3$, does not use the heat flow method and instead follows a different approach. There, the main work lies in establishing that the set of $K$-quasi\-conformal maps which admit a harmonic quasi-isometric extension is closed under pointwise convergence. The claim then follows from the existence result of \cite{LiTam93} for diffeomorphisms and the fact that every quasiconformal map of $\set{S}^2$ is a limit of uniformly quasiconformal diffeomorphisms. However, the analogue of the latter statement is not known for any higher-dimensional unit sphere and so we cannot use the same approach when $n\geq 4$. Nonetheless, there is some overlap with the methods used in \cite{Markovic14}. First and foremost, the good extension from Section \ref{sect:good} is a higher-dimensional analogue of the good extension from \cite{Markovic14}. Second, as already mentioned, Theorem \ref{thm:close} follows essentially from ideas in that paper.

What drives our proof behind the scenes is the quasiconformal Mos\-tow rigidity which holds in the hyperbolic space of dimension $ \ge 3$. More precisely, in order to prove that the good extension is almost harmonic at most points, see Proposition \ref{prop:cone}, we heavily use the fact that every quasiconformal map of $\set{S}^{n-1}$ with $n\geq 3$ is differentiable almost everywhere (and the derivative has maximal rank). It is known that Mostow rigidity fails for $n=2$ and consequently the existence proof in that case \cite{Markovic15} is very different from the ones in \cite{Markovic14} and here.

\section{Preliminaries}
We recall the following definitions.
\subsection{Quasi-isometries and quasiconformal maps} 

Let $F:X \to Y$ be a map between two metric spaces $(X,d_X)$ and $(Y,d_Y)$.

We say that $F$ is an $(L,A)$-quasi-isometry if there are constants $L>0$ and $A \ge 0$, such that
$$
\frac{1}{L}d_Y(F(x),F(y))-A \le \d_X(x,y) \le Ld_Y(F(x),F(y))+A,
$$
for every $x,y \in X$. 
An $(L,0)$-quasi-isometry is also called an $L$-Bi-Lipschitz map. 

We define the distortion function as
$$
\mathbf{K}(F)(x)=\limsup\limits_{t \to 0}  \,   \frac{ \max\limits_{d_{X}(x,y)=t} \d_Y(F(x),F(y)) }  {\min\limits_{\d_{X}(x,y)=t} \d_Y(F(x),F(y)) }.
$$
If $\mathbf{K}(F)(x) \le K$ on some  set $U \subset X$, we say that $F$ is locally $K$-quasiconformal on $U$. If $F$ is a global homeomorphism and $\mathbf{K}(F)(x) \le K$ for every $x \in X$, then we say that $F$ is $K$-quasiconformal (or $K$-qc for short).

We recall that every quasi-isometry $F:\set{H}^n\rightarrow \set{H}^n$ extends continuously to a quasiconformal map on $\del\set{H}^n\equiv\set{S}^{n-1}$. Two quasi-isometries $F,G$ have the same qc extension iff their distance $d_{\set{H}^n}(F(x),G(x))$ is uniformly bounded on $\set{H}^n$, see Proposition 1.6 in \cite{LiWang98}. 

\be{defn}
Let $a\in X$. We write $\mathbf{QC}_a(X)$ for the set of quasiconformal maps $F:X\rightarrow X$ which fix $a$, i.e.\ for which $F(a)=a$.
\e{defn}

For further background on quasi-isometries and qc maps, see \cite{LiWang98, TamWan98, Pansu}.

\subsection{Energy, tension field and harmonic maps} 
Let $(M,g)$, $(N,h)$ be Riemannian manifolds and let $F:M \to N$ be a $C^2$ map. The energy density of $F$ at a point $x \in M$ is defined as
$$
\mathbf{e}(F) = \frac{1}{2} |dF|^2
$$
where $|dF|^2$ is the cubed norm of the differential of $F$, taken with respect to the induced metric on the bundle $T^* M \times F^{-1} T N$. Equivalently,
$$
\mathbf{e}(F) = \frac{1}{2} \text{trace}_g F^* h
$$
and therefore in local coordinates
$$
\mathbf{e}(F) = \frac{1}{2} g^{ij} h_{\alpha\beta}\frac{\partial{F}^\alpha}{\partial x^i}\frac{\partial{F}^\beta}{\partial x^j},
$$
The tension field of $F$ is given by
$$
\tau(F)= \text{trace}_g \nabla dF,
$$
where $\nabla$ is the connection on the vector bundle $T^* M \times F^{-1} T N$  induced by the Levi-Civita connections on M and N.

$F$ is called \emph{harmonic} if $\tau(F)\equiv 0$. For background on harmonic maps see \cite{SchoenYau97, LinWang}.

\subsection{The heat equation}
Recall the heat equation with  initial map $\phi:\set{H}^n\rightarrow\set{H}^n$,
\beq
\label{eq:HE'}
\begin{aligned}
\tau(u)(x,t) &= \del_t u(x,t),\quad \textnormal{on } \set{H}^n\times [0,\it)\\
u(x,0)&=\phi(x),\quad \textnormal{on } \set{H}^n,
\end{aligned}
\eeq
A solution to the heat equation can be written in terms of the \emph{heat kernel} $H(x,y,t)$ as
\beq
\label{eq:HKdefn}
u(x,t)=\int\limits_{\set{H}^n} H(x,y,t) \phi(y) \d \lam(y),
\eeq
where $\d\lam$ is the volume measure for the hyperbolic metric. We quote a result which guarantees global in time existence and uniqueness of solutions to the heat equation for sufficiently nice initial maps $\phi$. It follows by combining Corollary 3.5 and Lemma 2.6 in \cite{Wang98}.

\be{prop}[Global in time existence and uniqueness]
\label{prop:HE}
Let $\phi:\set{H}^n\rightarrow\set{H}^n$ be a $C^2$-map with $\|\tau(\phi)\|\leq T$ for some $T>0$. Then, there exists a unique solution $u:\set{H}^n\times [0,\it)\rightarrow \set{H}^n$ to the heat equation \eqref{eq:HE} with initial map $\phi$.
\e{prop}

\section{Statement of two key results and proof of main result}

\subsection{Uniform convergence of the tension}
For $f\in\mathbf{QC}_a(\set{S}^{n-1})$, let $\curly{G}_a(f)$ be the good extension defined in Section \ref{sect:good}. Since $\curly{G}_a(f)$ is a $C^2$ map with uniformly bounded tension (see Definition \ref{defn:admissible} (ii)), Proposition \ref{prop:HE} implies that the heat equation with initial map $\phi\equiv\curly{G}_a(f)$ has a unique solution for all times, call it $u_{a,f}(x,t)$.

The following theorem is the first key result. It says that the tension field of $u_{a,f}(x,t)$ converges to zero, uniformly in space, as time goes to infinity.

\be{thm}
\label{thm:heat}
Let $n\geq 3$. For every $\eps>0$, there exists $t_0=t_0(K,\eps)$ such that for all $t\geq t_0$, we have 
\beq
\label{eq:heatclaim}
\|\tau (u_{a,f})(\cdot,t)\|<\eps,
\eeq
for every $a\in \set{S}^{n-1}$ and every $K$-qc map $f\in\mathbf{QC}_a(\set{S}^{n-1})$.
\e{thm}

The proof of Theorem \ref{thm:heat} is based on Hamilton's parabolic maximum principle and the new Sector Lemma \ref{lm:sector}. Here is a brief discussion of the ideas in the proof.

\be{itemize}
\item By Hamilton's parabolic maximum principle \cite{Hamilton75}, we have
\beqs
|\tau (u_{a,f})(x,t)|^2\leq \int\limits_{\set{H}^n} H(x,y,t) |\tau(\curly{G}_a(f))(y)|^2 \d \lam(y),
\eeqs
where $\d\lam$ denotes the volume measure for the hyperbolic metric.
\item We evaluate the integral in geodesic polar coordinates centered at $x$. We then use the diffusion of heat in hyperbolic space. Namely, we use that the heat kernel times the hyperbolic volume measure is effectively supported on a certain ``main annulus'' which travels to infinity as $t\rightarrow\it$. (We derive the main annulus in Appendix \ref{app:heat}, see Figure \ref{fig:Gaussian} for a picture.) 



\item Since the tension field of the good extension $\curly{G}_a$ is small at a random point, we expect
that $|\tau(\curly{G}_a(f))|$ becomes small on average on the main annulus. To prove the claim \eqref{eq:heatclaim}, though, we need this convergence to be uniform in $x$ (or, equivalently, uniform in $f$). This creates a problem since the heat dissipates in the hyperbolic space as the time goes to infinity.

\item The solution to this is to cover the main annulus by \emph{``good'' sectors} on which $|\tau(\curly{G}_a(f))|$ is small on average by the crucial Sector Lemma \ref{lm:sector}. Importantly, the good sectors have sizes which are bounded uniformly in $x$. As usual, uniformity is proved by appealing to the compactness of subsets of $K$-qc maps fixing certain points via Arzela-Ascoli.

\item To prove the Sector Lemma \ref{lm:sector}, it is helpful to work in a certain upper half space model of hyperbolic space. When choosing the upper half space model, other restrictions prevent us from also choosing which boundary point is mapped to infinity. Therefore, it is important for us that the good extensions at different boundary points are related via the \emph{partial conformal naturality} already mentioned in the introduction (see also Definition \ref{defn:pcn}).

\e{itemize}

\subsection{Every almost harmonic extension is close to the good extension}

The second key result is
\be{thm}
\label{thm:close}
Let $K\geq 1$ and $a\in\set{S}^{n-1}$. There exists $\eps_0=\eps_0(K)>0$ and $D_0=D_0(K,a)$ such that for all $K$-qc maps $f\in\mathbf{QC}_a(\set{S}^{n-1})$,
\beqs
\|d_{\set{H}^n}\l(\curly{G}_a(f),\psi\r)\|\leq D_0,
\eeqs
holds for all $C^2$ quasi-isometries $\psi:\set{H}^n\rightarrow\set{H}^n$ which extend $f$ and satisfy $\|\tau(\psi)\|<\eps_0$.
\e{thm}

The statement of Theorem \ref{thm:close} with $\eps_0=0$ was proved  in \cite{Markovic14}. The proof of Theorem \ref{thm:close} is word by word the same as the proof of the corresponding statement from \cite{Markovic14}, modulo their obvious generalization to higher dimensions and the observation that they provide sufficient ``wiggle room'' to allow for the presence of the $\eps_0$. There are two places where very minor changes have to  be made to the argument from \cite{Markovic14} and we will describe these below.

 Here is a very brief description of the proof of Theorem \ref{thm:close}. One uses the Green identity on a punctured ball to lower bound the maximum of
 $$\mathbf{d}^2\equiv d_{\set{H}^n}\l(\curly{G}_a(f),\psi\r)^2$$
  by the integral of its Laplacian times the Green's function. Then, one applies the usual lower bound on the Laplacian of the distance \cite{SchoenYau79, JagerKaul} to get a lower bound on this integral in terms of the maximum of $\mathbf{d}^2$ times a constant which depends on the radius of the ball. This constant can be made large by increasing the radius of the ball and one concludes that $\mathbf{d}^2$ is bounded.

\be{rmk}
In fact, with a little extra work the constant $D_0$ in Theorem \ref{thm:close} can be chosen independently of $a\in\set{S}^{n-1}$. To see this, one follows the same proof except that one replaces the compactness argument of Lemma 3.2 in \cite{Markovic14} with the slightly more elaborate one used in the proof of the Sector Lemma \ref{lm:sector} (i.e.\ essentially compactness of $\set{S}^{n-1}$ and continuity of the good extension $\G_a$ in $a$).
\e{rmk}


\subsection{Proof of main result}
\be{proof}[Proof of Theorem \ref{thm:conv} and Corollary \ref{cor:schoen}]
We assume Theorems \ref{thm:heat} and \ref{thm:close} hold. By conjugation with appropriate isometries and the conformal naturality of harmonic maps, it suffices to prove the claim for all $K$-qc $f\in\mathbf{QC}_a(\set{S}^{n-1})$ with $a\in \set{S}^{n-1}$ and $K\geq 1$ fixed.

By Proposition \ref{prop:admissible}, the good extension $\curly{G}_a(f)$ is admissible in the sense of Definition \ref{defn:admissible}. First, it is an $(L,A)$-quasi-isometry for some $L=L(K)$ and $A=A(K)$. Second, its tension is uniformly bounded, $\|\tau(\curly{G}_a(f))\|\leq T=T(K)$. From Hamilton's parabolic maximum principle (see \eqref{eq:tauestimate} below and recall that the integral of the heat kernel is normalized to one) we find that
$$\|\tau(u_{a,f})(\cdot,t)\|\leq T$$ for all $t\geq 0$. Since $u_{a,f}$ solves the heat equation, this implies that $\|d_{\set{H}^n}\l(\curly{G}_a(f),u_{a,f}(\cdot,t)\r)\|$ is bounded for every finite time $t$ (with a bound depending only on $t$ and $K$). By combining Theorems \ref{thm:heat} and \ref{thm:close}, the distance is also bounded for all $t\geq t_0(K)$. This proves the important intermediate result
\beq
\label{eq:almost}
\sup_{t>0}\|d_{\set{H}^n}\l(\curly{G}_a(f),u_{a,f}(\cdot,t)\r)\|\leq C
\eeq
for some constant $C=C(K)>0$. By a standard application of Cheng's Lemma \cite{Cheng} (see also \cite{Li,Wang98}), this gives a bound on the energy density of $u_{a,f}(x,t)$ which is uniform in $t$. This implies that $u_{a,f}(\cdot,t)$ and its derivative converges pointwise along some subsequence $t_i\rightarrow\it$ to a smooth map
\beqs
\curly{H}(f): \set{H}^n\rightarrow\set{H}^n,
\eeqs
which is harmonic by Theorem \ref{thm:heat}.
Recall that $\curly{G}_a(f)$ is an $(L,A)$-quasi-isometry and that by \eqref{eq:almost}, its distance to $\curly{H}(f)$ is bounded by $C=C(K)$. From this, we conclude that there exist $L_1=L_1(K)>0$ and $A_1=A_1(K)\geq0$ such that $\curly{H}(f)$ is an $(L_1,A_1)$-quasi-isometry, see e.g.\ \cite{Kapovich}. Finally, any two quasi-isometries which are at finite distance from each other have the same quasiconformal boundary map and therefore $\curly{H}(f)$ is an extension of $f$ as well. This proves Corollary \ref{cor:schoen}. Finally, by the uniqueness of the harmonic extension of a quasiconformal map \cite{LiWang98}, we can lift the subsequential convergence to convergence for all $t\rightarrow\it$. This proves Theorem \ref{thm:conv}.
\e{proof}

\section{Proof of Theorem \ref{thm:heat}}
The proof is based on Hamilton's parabolic maximum principle, the ballistic diffusion of heat discussed in hyperbolic space (see Appendix \ref{app:heat}) and the Sector Lemma \ref{lm:sector}. The first two facts are relatively well known. The Sector Lemma is at the heart of our argument, its proof is deferred to the next section.

\subsection{Hamilton's parabolic maximum principle and geodesic polar coordinates}
Fix $a\in\set{S}^{n-1}$. Recall the definition of the heat kernel $H(x,y,t)$ in \eqref{eq:HKdefn}. Since $\set{H}^n$ has negative curvature, it was observed by Hamilton \cite{Hamilton75} that
\beq
(\Del-\del_t)|\tau(u_{a,f})(x,t)|^2 \geq 0,
\eeq
for all $(x,t)\in \set{H}^n\times [0,\it)$. Hence, the parabolic maximum principle in the form of Theorem 3.1 in \cite{Wang98} implies
\beq
\begin{aligned}
\label{eq:tauestimate}
|\tau(u_{a,f})(x,t)|^2&\leq \int\limits_{\set{H}^n} H(x,y,t) |\tau(u_{a,f})(y,0)|^2 \d \lam(y)\\
&= \int\limits_{\set{H}^n} H(x,y,t) |\tau(\curly{G}_a(f))(y)|^2 \d \lam(y)
\end{aligned}
\eeq
for all $(x,t)\in \set{H}^n\times [0,\it)$.

The \emph{geodesic polar coordinates}  on $\set{H}^n$,  centered at $x \in \set{H}^n$, are given as follows. For $y \in \set{H}^n$, we write
\beqs
y=(\rho,\zeta)\qquad \text{with }\rho=d_{\set{H}^n}(x,y),\,\, \zeta\in\set{S}^{n-1},
\eeqs
where $\zeta$ is the unit vector at $x$ that is tangent to the geodesic ray  that starts at $x$ and contains $y$.  Using the standard identification between the unit tangent space at $x$ and the sphere $\set{S}^{n-1}$, we write $\zeta \in \set{S}^{n-1}$.

For a given $x\in\set{H}^n$, we will compute the integral on the right-hand side of \eqref{eq:tauestimate} in the geodesic polar coordinates centered at $x$. The volume element in the geodesic polar coordinates is 
\beqs
\sinh^{n-1}(\rho)\d\rho\,\d\zeta
\eeqs
 with $\d \zeta$ the Lebesgue measure on $\set{S}^{n-1}$. Since the heat kernel is a radial function (which we denote by $H(\rho,t)$), we have
\begin{align}
\nonumber
&\int\limits_{\set{H}^n} H(x,y,t) |\tau(\curly{G}_a(f))(y)|^2 \d \lam(y)\\
\label{eq:see}
&=\int\limits_0^\it   H(\rho,t) \sinh^{n-1}(\rho) \l(\,\int\limits_{\set{S}^{n-1}} |\tau(\curly{G}_a(f))(\rho,\zeta)|^2\d \zeta\r) \d \rho.
\end{align}
Next we will use the fact that heat travels approximately ballistically in the hyperbolic space (see Appendix \ref{app:heat}) to conclude that the radial integral in \eqref{eq:see} can be effectively restricted to a certain ($t$-dependent) ``main annulus''.

\subsection{Reduction to the main annulus}
Let $a\in\set{S}^{n-1}$ and let $f\in\mathbf{QC}_a(\set{S}^{n-1})$ be a $K$-qc map. For all $x\in\set{H}^n$, define the radial function
\beq
\label{eq:Phidefn}
\Phi(\rho):=\int\limits_{\set{S}^{n-1}} |\tau(\curly{G}_a(f))(\rho,\zeta)|^2\d \zeta,
\eeq
and recall that $\rho$ denotes hyperbolic distance from $x$. By Proposition \ref{prop:admissible}, $\{\curly{G}_a\}_a$ is an admissible family of extensions in the sense of Definition \ref{defn:admissible}. In particular, $\|\tau(\curly{G}_a(f))\|\leq T(K)\equiv T$. This implies that $\Phi$ is bounded and therefore it satisfies the assumption in Proposition \ref{prop:annulus} (iii).

We combine \eqref{eq:tauestimate}, \eqref{eq:see} and Proposition \ref{prop:annulus} (iii), which quantifies the extent to which the heat flow (more precisely the function $H(\rho,t)$ $\sinh^{n-1}(\rho)$) is concentrated on the main annulus. We express this as follows. Let $\eps'>0$. We find
\beq
\label{eq:see3}
|\tau(u_{a,f})(x,t)|^2 \leq \frac{C_n'}{\sqrt{t}} \int\limits_{R_{\mathrm{in}}}^{R_{\mathrm{out}}} \Phi(\rho) \d \rho  + \eps' ,
\eeq
for all $x\in\set{H}^n$ and for all $t\geq t_0=t_0(K,\eps')$. Here $C_n'$ is a universal (dimension dependent) constant and we introduced the inner and the outer radius of the main annulus
\beq
\label{eq:Rindefn}
\begin{aligned}
R_{\mathrm{in}}:=(n-1)t-l(\eps')\sqrt{t},\quad R_{\mathrm{out}}:=(n-1)t +l(\eps')\sqrt{t}.
\end{aligned}
\eeq
(We have $l(\eps')=\sqrt{8\log\l(\frac{C_n}{\eps'}\r)}$ where $C_n>0$ is another universal constant, but we will only need this formula in the final step of the proof.)

\subsection{Admissible sectors, good sectors and the Sector Lemma}
Recall that we write $(\rho,\zeta)$ for geodesic polar coordinates which are centered at $x\in\set{H}^n$. By a \emph{sector} we mean a set of the form
\beq
\label{eq:sectordefn}
S(x,\rho_{\mathrm{min}},r,\Om):=\setof{(\rho,\zeta) \in\set{H}^n}{\rho_{\mathrm{min}}\leq\rho\leq \rho_{\mathrm{min}}+r,\, \zeta\in \Omega}
\eeq
where $\rho_{\mathrm{min}},r>0$ and $\Omega\subset\set{S}^{n-1}$ (whenever we can, we write $S\equiv S(x,\rho_{\mathrm{min}},r,\Om)$).

  We will only consider the following class of \emph{admissible} sectors. Intuitively, a sector is admissible if the set $\Omega$ (of its ``angles'') has ``bounded geometry'', and if the diameter of $\Omega$ is comparable to $e^{-\rho_{\min}}$ (in particular,  note that admissibility is  independent of the choices of $x\in\set{H}^n$ and $r>0$ in the above definition \eqref{eq:sectordefn} of a sector). 

\be{defn}[Admissible sectors]
\label{defn:sadmissible}
Let $\al\geq 1$. We say that a sector  $S=S(x,\rho_{\mathrm{min}},r,\Om)$ if $\al$-admissible if  there exists a disk $D_{\n}\subset\set{S}^{n-1}$ of radius at least $\al^{-1} e^{-\rho_{\min}}$ and a disk $D_\out\subset\set{S}^{n-1}$ of radius at most $\al e^{-\rho_{\min}}$ (both in $\set{S}^{n-1}$ distance) such that
\beq
D_\n \subset \Om\subset D_\out.
\eeq
In this case, it will be convenient to call $\Om$ an $(\al,\rho_{\min})$-admissible subset of $\set{S}^{n-1}$.
\e{defn}

The only example of a 1-admissible sector is when the corresponding set $\Om\subset\set{S}^{n-1}$ is a disk of radius $e^{-\rho_{\min}}$ in $\set{S}^{n-1}$ distance (i.e.\ a small spherical cap).  Generalizing this example to $\al$-admissible sectors will give us extra flexibility when we apply the Sector Lemma in the next section (it is easier to ``stack'' admissible sectors if the $\Om$ do not have to be exactly disks). 

The Sector Lemma \ref{lm:sector} below is at the heart of our proof. It says that for a given $\al\geq 1$, and when   $\rho_{\mathrm{min}}$ is sufficiently large, there exists $r_1>0$ such that  every $\al$-admissible  sector $S=S(x,\rho_{\mathrm{min}},r_1,\Om)$ is ``good'' in the sense that  the  tension field (of the good extension) is \emph{small on average} on $S$. We first define a good sector.

\be{defn}[Good sectors]
\label{defn:sgood} Let $\de>0$, $a\in\set{S}^{n-1}$ and let $f\in\mathbf{QC}_a(\set{S}^{n-1})$ be a $K$-qc map. We say that a sector $S$ (as defined by (\ref{eq:sectordefn})) 
is $\delta$-\emph{good} (or just \emph{good} if $\delta$ is understood), if
\beq
\int\limits_{S} |\tau(\G_a(f))(\rho,\zeta)|^2 \d\rho\d\zeta
< \de \int\limits_{S} 1\, \d\rho\d\zeta.
\eeq
\e{defn}
Next, we state the Sector Lemma. But before we do that, recall that the notion of admissibility of a sector $S(x,\rho_{\mathrm{min}},r,\Om)$ does not depend on the choice of $r>0$ (it also does not depend on the choice of $x\in\set{H}^n$ but we will not use this here). In other words, given $r,r'>0$, we have that the sector $S(x,\rho_{\mathrm{min}},r,\Om)$ is $\al$-admissible if and only if the sector  $S(x,\rho_{\mathrm{min}},r',\Om)$ is $\al$-admissible. By $\{S(x,\rho_{\mathrm{min}},r,\Om)\}_{r}$ we denote the family of sectors when $r$ varies over $(0,\it)$, and we say that this family of sectors is $\al$-admissible if all of the sectors (or equivalently one of them) are $\al$-admissible.

\be{lm}[Sector Lemma]
\label{lm:sector} Let $\al, K\geq 1$ and $\de>0$. There exist constants $r_0=r_0(K,\al,\de)>1$ and $\rho_0=\rho_0(K,\al,\de)>0$ such that for all $x\in\set{H}^n$, $a\in\set{S}^{n-1}$ and for all $K$-qc maps $f\in\mathbf{QC}_a(\set{S}^{n-1})$ the following holds.  Every  $\al$-admissible  family $\{S(x,\rho_{\mathrm{min}},r,\Om)\}_{r}$ which satisfies $\rho_{\mathrm{min}}>\rho_0$,  contains a $\delta$-good sector $S(x,\rho_{\mathrm{min}},r_1,\Om)$, where $1\leq r_1 \leq r_0$.
\e{lm}

The Sector Lemma will be proved later, see Section \ref{sect:sp}. For now, we continue with the proof of Theorem \ref{thm:heat}. 

Before we go on with this, we remark why the factor $e^{-\rho_{\min}}$ appears in Definition \ref{defn:sadmissible} of an admissibile sector.

\be{rmk}
The factor $e^{-\rho_{\min}}$ in Definition \ref{defn:sadmissible} is important in the proof of the Sector Lemma \ref{lm:sector}, which is given in Section 6. The proof seeks to get a contradiction to the existence of a sequence of ``bad'' sectors which will have to ``run off'' to the boundary of hyperbolic space (i.e.\ ${\rho_{\min}}\rightarrow\it$). Going to appropriate boundary coordinates leads to the angular variable being rescaled by a factor proportional to $e^{\rho_{\min}}$. To get a contradiction, we need the sequence of bad sectors to yield a limiting set which has diameter of order one and this is only possible if we initially scale down the angular variable by $e^{-\rho_{\min}}$.
\e{rmk}

\subsection{Covering the main annulus with good sectors}
Recall that $\eps'$ is a fixed positive quantity, which we will eventually let  go to zero. Our goal is to estimate the right hand side in \eqref{eq:see3}, i.e.\ the average of $|\tau(\G_a(f))|^2$ over the main annulus, by the small quantity $\eps'$.

We will achieve this by covering the main annulus with $\eps'$-good sectors (i.e.\ sectors on which $|\tau(\curly{G}_a(f))|$ is small on average, see Definition \ref{defn:sgood}). Such $\eps'$-good sectors exist by the Sector Lemma \ref{lm:sector} (in the following we will just speak of ``good'' sectors, $\eps'$ is understood). 

We cover the main annulus with good sectors in two steps. In step 1, we cover the main annulus by cylinders (in geodesic polar coordinates) which do not overlap too much. In step 2, we partition each cylinder (up to a small region near its top) into good sectors. This partitioning is most conveniently performed when the cylinders are mapped to Euclidean cuboids and our notion of an admissible sector is flexible enough to allow for this.

\subsection{Step 1: Covering the main annulus by cylinders}
We first note that one can of course cover the sphere efficiently by small disks. 

 \be{prop}
\label{prop:besi}
There exists a universal constant $\beta_n>0$ (the universal constant from the Besicovitch covering theorem in $\set{R}^n$) such that the following holds. For every $R>0$, there is a finite covering $\{D_i\}_{1\leq i\leq M}$ ($M$ is some finite integer) of $\set{S}^{n-1}$ by disks $D_i\subset\set{S}^{n-1}$ of radius $e^{-R}/2$ (in $\set{S}^{n-1}$ distance) such that every point of $\set{S}^{n-1}$ is contained in at most $\beta_n$ of the disks.
\e{prop}

\be{proof}
We cover $\set{S}^{n-1}$ by taking a disk of radius $e^{-R}/2$ (in $\set{S}^{n-1}$-distance) around every point. By compactness, we can pass to a finite subcover. By the Besicovitch covering theorem, there exists a universal constant $\beta_n$ and yet another finite subcover, call it $\{D_i\}_{1\leq i\leq M}$ such that 
every point on the sphere is contained in at most $\beta_n$ of the $D_i$. (Formally, one first takes the finite subcover which exists by compactness and extends the disks to get a finite covering of $\set{S}^{n-1}$ by $n$ dimensional balls, with centers on $\set{S}^{n-1}$. One then applies the Besicovitch covering theorem in $\set{R}^n$ to these balls and replaces them by the corresponding disks again.)
\e{proof}

We recall that the main annulus from Proposition \ref{prop:annulus} is of the form (in geodesic polar coordinates)
\beq
[R_{\mathrm{in}},R_{\mathrm{out}}]\times \set{S}^{n-1}
\eeq
where the inner and the outer radius are given by \eqref{eq:Rindefn}.

We now apply Proposition \ref{prop:besi} with the choice $R=R_{\n}$. This yields a covering of $\set{S}^{n-1}$ by disks 
\beqs
\{D_i\}_{1\leq i\leq M}
\eeqs
of radius $e^{-R_{\mathrm{in}}}/2$ (in $\set{S}^{n-1}$ distance) such that every point in $\set{S}^{n-1}$ is contained in at most $\beta_n$ of the disks. For every $1\leq i\leq M$, we define the cylinder (in geodesic polar coordinates)
\beq
\label{eq:Cidefn}
\curly{C}_i:=[R_{\mathrm{in}},R_{\mathrm{out}}] \times D_i.
\eeq
Each cylinder $\curly{C}_i$ covers the portion of the main annulus which has ``angle'' $\zeta\in D_i$. Notice that each point in the main annulus lies in at most $\beta_n$ of the cylinders $\curly{C}_i$ (because the same holds true for the disks $D_i$).


\subsection{Step 2: Partitioning the cylinders into good sectors}

In step 2, we shall partition each cylinder $\curly{C}_i$ into good sectors (excluding a small region near the top of the cylinder). Good sectors exist by the Sector Lemma \ref{lm:sector}. 

The idea is \emph{to apply the Sector Lemma \ref{lm:sector} iteratively}. That is, starting at $\rho_{\min}=R_{\n}$, we \emph{stack good sectors on top of each other} until we (almost) reach $R_{\out}$.  The process is as follows. Once we have added a good sector to the partition of $\curly{C}_i$, we then partition the top of this sector into admissible domains (the last sentence in Definition \ref{defn:sadmissible}), and then erect a good sector above each admissible domains. We then partition the top of each new sector and so on. We stop adding good sectors when the total height of a stack gets close to $R_{\out}$, so as not to overshoot.

\begin{figure}[t]
\be{center}
\includegraphics[height=9cm]{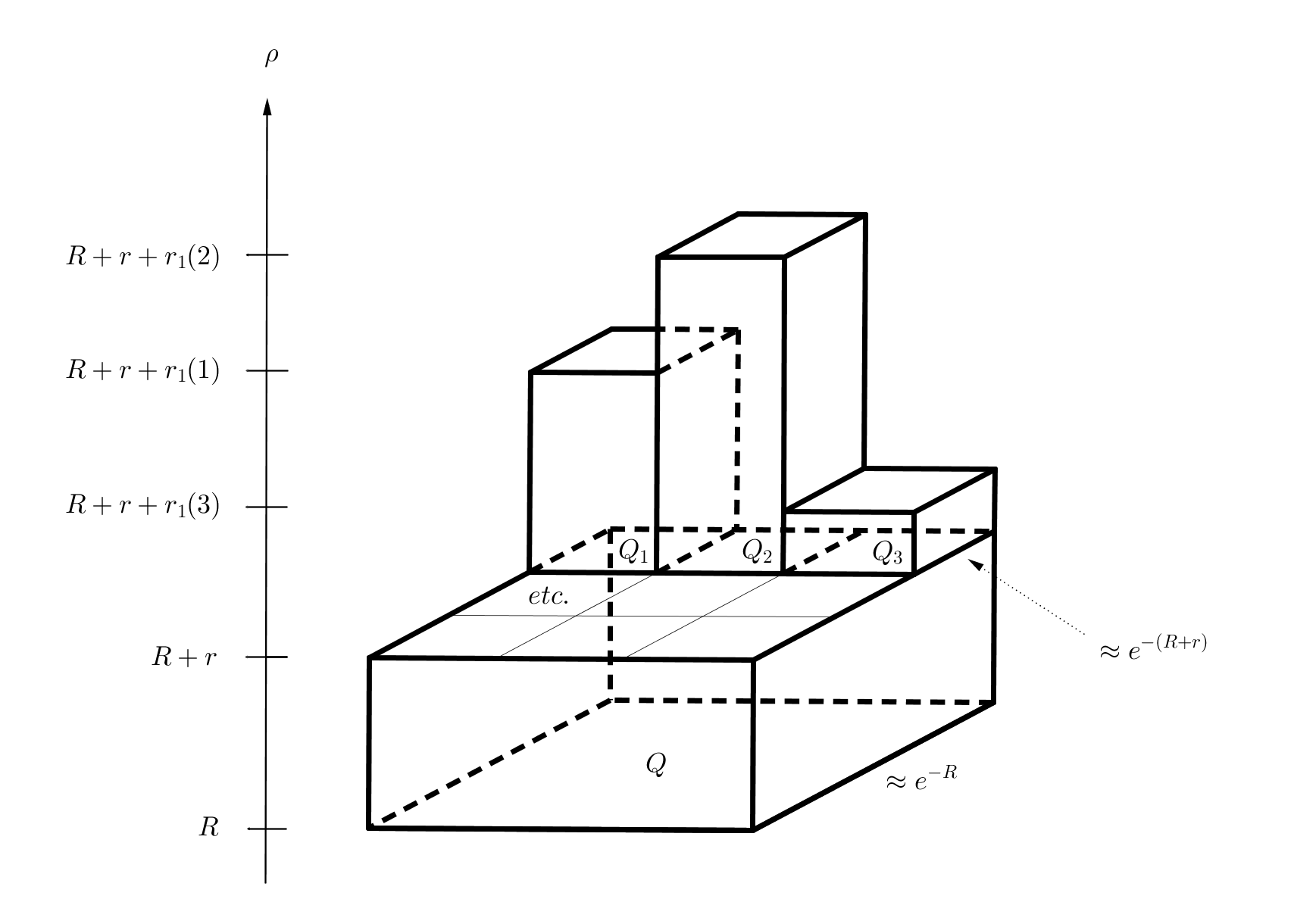}
\caption{This picture shows one step in our inductive partitioning of the cylinder $\curly{C}_i$ into good sectors. We partition the top face of a given cuboid $[R,R+r]\times Q$ and erect a new cuboid on top of each subface $Q_j$. The new cuboid has the height $r_1(j)$ determined by the Sector Lemma. To obtain the new good \emph{sectors}, each Euclidean cube $Q_j$ is mapped to some $B(Q_j)\subset\set{S}^{n-1}$ by a uniformly Bi-Lipschitz map $B$. Notice that each $B(Q_j)$ will be \emph{admissible} (see Definition \ref{defn:sadmissible}) in the right way, because it is the Bi-Lipschitz image of a cube $Q_j$ with the correct sidelength $\approx e^{-(R+r)}$ (here $\approx$ means equality up to a factor of two).
}
\label{fig:cubes}
\e{center}
\end{figure}

It is important that each new admissible sector is $\al$-admissible, where $\al$ is some universal constant. Therefore at each inductive step we are required to partition an 
admissible domain into  admissible domains of an appropriate (smaller) diameter so that the new domains have uniformly bounded geometry . This is easily achievable when the domain we want to partition is a Euclidean cube, see Figure \ref{fig:cubes}. 

Thus, it is most convenient to stack the good sectors on top of each other when their base (originally a subset of the sphere) can be viewed as a Euclidean cube in $\set{R}^{n-1}$.  We achieve this by mapping each cylinder $\curly{C}_i=[R_{\mathrm{in}},R_{\mathrm{out}}] \times D_i$ using a  uniformly Bi-Lipschitz map onto the Euclidean cuboid 
$[R_{\mathrm{in}},R_{\mathrm{out}}] \times E(R_\mathrm{in})$, where $E(R_\mathrm{in})\subset\set{R}^{n-1}$ is the Euclidean cube of diameter 
$e^{-R_\mathrm{in}}$ and centred at the origin. We then perform the partition in the cuboid model and return it back to  $\curly{C}_i$ with the Bi-Lipschitz map.

The upshot is (recall that $\eps'>0$ is fixed)
\be{lm}
\label{lm:partition}
Let $1\leq i\leq M$ and $x\in\set{H}^n$. There exists $t_0=t_0(K,\eps')$ such that for all $t\geq t_0$, there exists a finite collection $\{S^{(i)}_j\}_{1\leq j\leq J}$ of disjoint sectors (a sector is a set of the form \eqref{eq:sectordefn})
that is contained in $\curly{C}_i$ and almost covers $\curly{C}_i$, i.e.
\beq
\label{eq:lm1}
\int\limits_{\curly{C}_i\setminus \bigsqcup\limits_{1\leq j\leq J}S_j^{(i)}} 1\, \d\rho\d\zeta < r_0(K,\eps') \int\limits_{D_i} 1\, \d\zeta
\eeq
where $r_0$ is defined by the Sector Lemma \ref{lm:sector}. Moreover, the sectors are $\eps'$-good in the sense of Definition \ref{defn:sgood}, i.e.
\beq
\label{eq:lm2}
\int\limits_{S^{(i)}_j} |\tau(\G_a)(f)(\rho,\zeta)|^2
\d\rho\d\zeta< \eps' \int\limits_{S^{(i)}_j} 1\,
\d\rho\d\zeta.
\eeq
\e{lm}

\be{rmk}
In fact, we will see in the proof below that the sectors $\{S^{(i)}_j\}_{1\leq j\leq J}$ are $\al$-admissible for some universal constant $ \al>1$, and this is why we drop the dependence of $r_0$ on the constant $\al$
characterizing the admissibility of the good sector from Sector Lemma \ref{lm:sector}.
\e{rmk}

We now give  formal proofs following the ideas sketched above. 

\be{proof}[Proof of Lemma \ref{lm:partition}]
Fix $1\leq i\leq M$ and $x\in\set{H}^n$. For simplicity we let $\curly{C}_i=\curly{C}$ and $D_i=D$. Let $E(R_\mathrm{in})\subset\set{R}^{n-1}$ denote the Euclidean cube of sidelength $e^{-R_\mathrm{in}}$ centered at the origin of $\set{R}^{n-1}$. There exists a Bi-Lipschitz map 
\beqs
B:E(R_\mathrm{in})\rightarrow D
\eeqs
with a Bi-Lipschitz constant bounded by a universal (dimension dependent) constant $L_0>1$. (This holds because the disk and the cube both have diameters which are proportional to $e^{-R_{\n}}$ up to a universal dimension dependent factor. Note also that this diameter is small, so that the disk $D\subset\set{S}^{n-1}$ is almost flat.) 

We now define the partition of the cylinder $\curly{C}$ into good sectors by apply the Sector Lemma \ref{lm:sector} inductively. In every application of the Sector Lemma, we shall choose $\de=\eps'$ (which was fixed before) and $\al=\sqrt{n}L_0>1$. Since $L_0$ is a universal constant, the quantities $r_0(K,L_0,\eps')>1$ and $\rho_0(K,L_0,\eps')>0$ provided by the Sector Lemma only depend on $K,\eps'$. By choosing $t\geq t_0$ with $t_0=t_0(K,\eps')$ sufficiently large, we can ensure that 
$$
R_{\mathrm{in}}:=(n-1)t-l(\eps')\sqrt{t}\geq \rho_0(K,\eps')
$$ 
holds for all $t\geq t_0$ (this is important because we want to choose $\rho_{\mathrm{min}}=R_{\mathrm{in}}$ next).

The \emph{inductive base case} is the following.  We apply the Sector Lemma \ref{lm:sector} with $\rho_{\min}=R_\mathrm{in}$ and $\Om=D$, which we note is $(1,R_{\n})$-admissible in the sense of Definition \ref{defn:sadmissible}  because $D\subset\set{S}^{n-1}$ is a disk of radius $e^{-R_{\n}}$. The Sector Lemma then says that the  sector 
$S_1 \equiv S(x,R_{\n},r_1,D)$ is  $\eps'$-good for some  $1\leq r_1 \leq r_0(K,\eps')$ (here we use that  $\Om=D$  is in particular $(\sqrt{n}L_0,R_{\n})$-admissible in the sense of Definition \ref{defn:sadmissible} so that we can apply the Sector Lemma with the corresponding $r_0=r_0(K,\epsilon')$ defined above). The sector $S_1$ is $\eps'$-good in the sense of Definition \ref{defn:sgood} and thus it satisfies \eqref{eq:lm2}. Equivalently, this sector can be written as
\beqs
S_1=[R_{\n},R_{\n}+r_1]\times D
\eeqs
which is  the first layer of the required partition of the cylinder $\curly{C}$. But most importantly from the point of view of our induction process, we note that one can also write
\beqs
S_1=[R_{\n},R_{\n}+r_1]\times B \big( E(R_\mathrm{in}) \big),
\eeqs
where we recall that $E(R_\mathrm{in})$ is the cube whose side length is $e^{-R_\mathrm{in}}$.

For what follows, the reader may find it helpful to consider Figure \ref{fig:cubes}. The \emph{inductive hypothesis} is the following. Suppose that an $\eps'$-good sector $S \equiv S(x,R,r,\Omega)$ is  included in the partition of the cylinder $\curly{C}$. Here we assume that $R_{\n} \leq R$ and  $1< r \leq R_{\out}-R$, and that 
\beqs
S=[R,R+r] \times B(Q) ,
\eeqs
where $Q \subset E(R_{\n})$ is a cube of side length between $e^{-R}$ and $2e^{-R}$ (note that it follows from these induction hypotheses that such a sector is $\sqrt{n}L_0$-admissible since $B$ is $L_0$ Bi-Lipschitz and since the sidelength of $Q$  belongs to the interval  $[e^{-R},2e^{-R}]$).

The \emph{inductive step} is as follows. If $R+r>R_{\out}-r_0(K,\eps')$ we stop. If not, we partition $Q$ into Euclidean cubes $Q_1,Q_2,...,Q_N$,  which all have the same sidelength that lives in the interval  $[e^{-(R+r)},2e^{-(R+r)} ]$ (it is elementary to see that such a partition of $Q$ always exists when $r \geq 1$).  

We include the following new sectors into the partition of $\curly{C}$. For $1\leq j \leq N$,  we let $S_j \equiv S(x,R+r,r_1(j),B(Q_j))$, where $1 \leq r_1(j) \leq r_0(K,\epsilon')$ is given by the Sector Lemma so that $S_j$ is an $\epsilon'$-good sector. Note that $S_j$  is $\sqrt{n}L_0$-admissible since $B$ is $L_0$ Bi-Lipschitz and since the sidelength of $Q_j$  belongs to the interval  $[e^{-(R+r)},2e^{-(R+r)}]$, and so we  can apply the Sector Lemma with the corresponding $r_0=r_0(K,\epsilon')$ defined above.

The new sectors $S_j$ satisfy the inductive hypothesis and we continue the induction until we have that $R+r>R_{\out}-r_0(K,\eps')$ for a sector $S \equiv S(x,R,r,\Omega)$ that is in the partition. Since each time when we add a new sector we increase the height by at least 1 (recall that $r_1$ from the Sector Lemma is at least 1), we will stop adding new sectors after finitely many steps. Since the sectors that partition $\curly{C}$ were all chosen to be $\eps'$-good in the sense of Definition \ref{defn:sgood} the relation  \eqref{eq:lm2} is immediate.
\e{proof}

 \subsection{Conclusion}
 We will now use the covering of the main annulus by good sectors to estimate \eqref{eq:see3}, i.e.\ the integral of the tension field of the good extension over the main annulus. This is the last step in proving Theorem \ref{thm:heat}.
 \be{proof}[Proof of Theorem \ref{thm:heat}]
Recall \eqref{eq:see3}
\beq
\label{eq:sec33}
|\tau(u_{a,f})(x,t)|^2 \leq \frac{C_n'}{\sqrt{t}} \int\limits_{R_\mathrm{in}}^{R_\mathrm{out}} \Phi(\rho) \d \rho  + \eps'.
\eeq
where we used the notation defined in \eqref{eq:Rindefn}. Recall that $\Phi(\rho)$ is defined in \eqref{eq:Phidefn} as the spherical average of the function $|\tau(\G_a(f))|^2$. Since this function is non-negative, we can estimate the integral over the main annulus by the integral over its covering $\cup_{i=1}^M\curly{C}_i$, where the cylinders $\curly{C}_i$ are defined in \eqref{eq:Cidefn}. This gives 
 \beq
 \label{eq:rewrite}
 \frac{1}{\sqrt{t}}\int\limits_{R_\mathrm{in}}^{R_\mathrm{out}} \Phi(\rho) \d \rho
 \leq \frac{1}{\sqrt{t}}\sum_{i=1}^M \int\limits_{\curly{C}_i} |\tau(\G_a(f))(\rho,\zeta)|^2 \d\rho\d\zeta. 
 \eeq
We now estimate this using Lemma \ref{lm:partition}. We first apply \eqref{eq:lm1}, i.e.\ we estimate the integral over each cylinder $\curly{C}_i$ by the integral over the finite disjoint union of good sectors $\bigsqcup\limits_{1\leq j\leq J} S^{(i)}_j$, up to a small region on which we bound the tension field by $T$. Then, we use that the sectors are $\eps'$-good, i.e.\ the average of the tension field is small on them, see \eqref{eq:lm2}. We get
\beq
\label{eq:paren}
\begin{aligned}
&\frac{1}{\sqrt{t}}\sum_{i=1}^M \int\limits_{\curly{C}_i} |\tau(\G_a(f))(\rho,\zeta)|^2 \d\rho\d\zeta\\
&\leq \frac{1}{\sqrt{t}}\sum_{i=1}^M \l(\,\int\limits_{\bigsqcup\limits_{1\leq j\leq J} S^{(i)}_j}|\tau(\G_a(f))(\rho,\zeta)|^2 \d\rho\d\zeta +r_0 T^2 \int\limits_{D_i}  1\, \d\zeta \r)\\
&\leq \frac{1}{\sqrt{t}}\sum_{i=1}^M \l(\,\eps' \sum_{1\leq j\leq J}\int\limits_{S^{(i)}_j} 1\, \d\rho\d\zeta + r_0 T^2 \int\limits_{D_i}1\, \d\zeta \r).  
\end{aligned} 
\eeq
Recall from Proposition \ref{prop:besi} that the disks do not overlap too much: For every point in $\set{S}^{n-1}$ is contained in at most $\beta_n$ of the disks $D_i$ (and $\beta_n$ is a universal constant). First, this gives
\beqs
\frac{r_0 T^2 }{\sqrt{t}}\sum_{i=1}^M \int\limits_{D_i} 1\, \d\zeta \leq \frac{r_0 T^2 }{\sqrt{t}} \beta_n |\set{S}^{n-1}|,
\eeqs
where $|\cdot|$ denotes the Lebesgue measure. Moreover, we recall that the good sectors are contained in the cylinder
\beqs
\bigsqcup\limits_{1\leq j\leq J} S^{(i)}_j\subset \curly{C}_i\equiv [R_{\mathrm{in}},R_{\mathrm{out}}]\times D_i
\eeqs
 and then we use that every point in $\set{S}^{n-1}$ is also contained in at most $\beta_n$ of the cylinders $\curly{C}_i$ to get
\beqs
\begin{aligned}
&\frac{\eps'}{\sqrt{t}}\sum_{i=1}^M \sum_{1\leq j\leq J}\int\limits_{S^{(i)}_j} 1\, \d\rho\d\zeta
=\frac{\eps'}{\sqrt{t}}\sum_{i=1}^M \int\limits_{\bigsqcup\limits_{1\leq j\leq J} S^{(i)}_j} 1\, \d\rho\d\zeta\\
\leq &\frac{\eps'}{\sqrt{t}}\sum_{i=1}^M \int\limits_{\curly{C}_i} 1\, \d\rho\d\zeta
\leq  \frac{\eps'}{\sqrt{t}} \beta_n|\set{S}^{n-1}| \int\limits_{R_{\mathrm{in}}}^{R_{\mathrm{out}}} 1\, \d\rho\\
= &2\eps' l(\eps') \beta_n|\set{S}^{n-1}|.
\end{aligned}
\eeqs
In the last step, we used that $R_\mathrm{out}-R_\mathrm{in}=2l(\eps')\sqrt{t}$, see their definition \ref{eq:Rindefn}. Recall that $l(\eps')=\sqrt{8\log\l(\frac{C_n}{\eps'}\r)}$ where $C_n>0$ is a universal constant. Combining \eqref{eq:sec33}-\eqref{eq:paren} and the estimates following them, we have shown that
\beq
\label{eq:UB}
|\tau(u_{a,f})(x,t)|^2 \leq 2\eps' l(\eps') C_n' \beta_n|\set{S}^{n-1}|+ \frac{C_n' r_0 T^2 \beta_n |\set{S}^{n-1}|}{\sqrt{t}} +\eps'
\eeq
We can now take the $\sup_{x\in\set{H}^n}$ on both sides (the right hand side no longer depends on $x$). The second term on the right hand side can be made less than $\eps'$ for all $t\geq t_0$ and $t_0=t_0(K,\eps')$ sufficiently large (recall that $r_0=r_0(K,\eps')$, $T=T(K)$ and $C_n',\beta_n$ are universal constants). 

Finally, observe that because $l(\eps')=\sqrt{8\log\l(\frac{C_n}{\eps'}\r)}$ with $C_n$ a universal constant, the first term on the right hand side in \eqref{eq:UB} vanishes as $\eps'\rightarrow0$. This proves Theorem \ref{thm:heat}.
\e{proof}

\section{The good extension of a quasiconformal map}
In this section, we discuss the good extension and its properties in some detail. First, we define the good extension $\G_\it(f)$ as in \cite{Markovic14} for quasiconformal boundary maps $f:\set{R}^{n-1}\rightarrow \set{R}^{n-1}$ which fix $\it$ (in the upper half-space model of hyperbolic space). We observe some of its important properties, in particular that $\G_\it$ is partially conformally natural with respect to isometries which fix $\it$, see Proposition \ref{prop:GE} (i). Then we extend the definition of the good extension to quasiconformal boundary maps $f:\set{S}^{n-1}\rightarrow\set{S}^{n-1}$ which fix an arbitrary point $a\in \set{S}^{n-1}$. Importantly, the resulting family of good extensions $\{\G_a\}_a$ satisfies partial conformal naturality (see Definition \ref{defn:pcn}) and it is admissible in the sense of Definition \ref{defn:admissible} (in particular it is continuous in $a$).

\label{sect:good}
\subsection{Preliminaries}
First we work in the upper half space model of hyperbolic space in Euclidean coordinates
\beqs
\mathbf{H}^n =\setof{(x,s)}{x\in\set{R}^{n-1},\, s>0}.
\eeqs
We identify $\del\mathbf{H}^n\equiv \ol{\set{R}^{n-1}}$ in the natural way. Recall that we write $\mathbf{QC}_\it(\ol{\set{R}^{n-1}})$ for the set of quasiconformal maps $\ol{\set{R}^{n-1}}\rightarrow \ol{\set{R}^{n-1}}$ which fix $\it$.
By the quasiconformal \emph{Mostow rigidity}, every such $f$ is differentiable almost everywhere (with the derivative of maximal rank). The energy density of $f\in\mathbf{QC}_\it(\ol{\set{R}^{n-1}})$ with respect to the Euclidean metric is then defined almost everywhere and reads
\beqs
 \mathbf{e}(f)(x)=\sum_{i,j=1}^{n-1} \sum_{\al,\beta=1}^{n-1} \frac{\del f_i}{\del x_j} \frac{\del f_\al}{\del x_\beta},
\eeqs
where we wrote $f=(f_1,\ldots, f_{n-1})$. We now define the good extension of all maps $f\in\mathbf{QC}_\it(\ol{\set{R}^{n-1}})$. We use the higher-dimensional analogue of the definition in \cite{Markovic14}, compare also \cite{BeurlingAhlfors}, \cite{KovalevOnninen}.

\be{defn}
\label{defn:gtilde}
For $f\in \mathbf{QC}_\it(\ol{\set{R}^{n-1}})$, define its good extension $\curly{G}(f):\mathbf{H}^n\rightarrow \mathbf{H}^n$ by
\beq
\begin{aligned}
\label{eq:gtildedefn}
\G_\it(f)(x,s):=&\l(\,\int\limits_{\set{R}^{n-1}} f(x+sy) \phi(y) \d y\r.,\\ 
&\quad\,\l.\frac{s}{\sqrt{n-1}}\sqrt{\int\limits_{\set{R}^{n-1}} \mathbf{e}(f)(x+sy) \phi(y)\d y}\r),
\end{aligned}
\eeq
where $\phi$ is the standard Gaussian
\beq
\label{eq:Gaussian}
\phi(y):=(2\pi)^{\frac{1-n}{2}}e^{-\frac{|y|^2}{2}}.
\eeq
\e{defn}

We write $\mathrm{Isom}_\it(\mathbf{H}^n)$ for the subset of isometries which fix $\it\in\ol{\set{R}^{n-1}}$. Note that
\beq
\label{eq:Isominfty}
\begin{aligned}
\mathrm{Isom}_\it(\mathbf{H}^n)=\{&(x,s) \mapsto(a O(x)+b,a s)\, :\\  
&a>0,\, b\in \set{R}^{n-1},\, O\in SO(n-1)\}
\end{aligned}
\eeq

As in \cite{Markovic14}, the good extension $\G_\it$ has the following properties. Unlike in \cite{Markovic14}, the partial conformal naturality from (ii) will be very important for us.

\be{prop}
\label{prop:GE}
For all $f\in \mathbf{QC}_\it(\ol{\set{R}^{n-1}})$, the integrals in \eqref{eq:gtildedefn} are well-defined and $\G_\it(f)\in C^\it(\mathbf{H}^n)$.
\be{enumerate}[label=(\roman*)]
\item $\G_\it$ is partially conformally natural under isometries fixing infinity, i.e.\
\beqs
\G_\it(I \circ f\circ J)= I\circ\G_\it(f)\circ J
\eeqs
for any $I,J\in \mathrm{Isom}_\it(\mathbf{H}^n)$.
\item Let $\curly{L}(\set{R}^{n-1})$ denote the set of invertible, orientation preserving linear maps from $\set{R}^{n-1}$ to itself. For every $L\in\curly{L}(\set{R}^{n-1})$, $\G_\it(L):\mathbf{H}^n\rightarrow \mathbf{H}^n$ is harmonic and satifies
\beq
\label{eq:linproperties}
 \mathbf{e}(\G_\it(L))(x,s)>1,\quad \mathbf{K}(\G_\it(L))(x,s)=\mathbf{K}(L)(x),
 \eeq
 for all $(x,s)\in \mathbf{H}^n$.
\e{enumerate}
\e{prop}

\be{proof}

The fact that the good extension is well defined and smooth follows by analogous arguments as in \cite{Markovic14}.

Statement (i) can be checked explicitly from \eqref{eq:gtildedefn} and \eqref{eq:Isominfty} as well as normalization and rotational invariance of the Gaussian. 

For statement (ii), we use a result of \cite{LiTam93} (see also \cite{TamWan98}), namely that every $L\in\curly{L}(\set{R}^{n-1})$ has a harmonic quasi-isometric extension which is given by
 \beqs
 \l(L(x),\sqrt{\frac{\mathbf{e}(L)}{n-1}}s\r).
 \eeqs
 It is elementary to check that $\G_\it(L)(x,s)$ defined by \eqref{eq:gtildedefn} takes precisely this form when $f\equiv L$ is linear. Therefore, $\G_\it(L)$ is harmonic. The properties \eqref{eq:linproperties} follow as in \cite{Markovic14}.
\e{proof}

The next statement is a slight (and straightforward) strengthening of Theorem 3.1 in \cite{Markovic14} to cones (extending into $\mathbf{H}^n$ starting from a tip in $\ol{\set{R}^{n-1}}$). In particular, it shows that eventually (as one moves towards the boundary of hyperbolic space) $\G_\it(f)$ is almost harmonic for any $f\in \mathbf{QC}_\it(\ol{\set{R}^{n-1}})$. 

\be{prop}
\label{prop:cone}
For $\eps>0$ and a $K$-qc map $f\in\mathbf{QC}_\it(\set{S}^{n-1})$, define the ``good set'' by
\beqs
\begin{aligned}
X_f(\eps):=\l\{(x,s)\in\mathbf{H}^n\, :\, \r. &\mathbf{e}(\G_\it(f))(x,s)>1,\, \mathbf{K}(\G_\it(f))(x,s)<2K,\\
&\l. |\tau(\G_\it(f))(x,s)|<\eps\r\}.
\end{aligned}
\eeqs
Then, for almost every $x\in\set{R}^{n-1}$,
\beqs
\lim_{s\rightarrow0}\,\left(\min\limits_{x':|x-x'|\leq s} \ind_{X_f(\eps)}(x',s) \right)=1,
\eeqs
where $\ind$ denotes the characteristic function of a set.
\e{prop}
This proposition says that for almost every $x\in\set{R}^{n-1}$, the geodesic  (together with the cone around it) starting at $\infty$ and ending at $x$ will eventually be contained in the good set $X_f$.
More precisely,  for almost every $x$ there exists a vertical (in the Euclidean sense)  geodesic ray ending at $x$ that together with the equidistant cone around it is contained in the good set $X_f$.

We will use this proposition on two occasions: (a) At the end of the proof of the Sector Lemma \ref{lm:sector}, we use that the tension field becomes small on the whole cone. (b) When following the arguments in \cite{Markovic14} to prove Theorem \ref{thm:close} (there we do not need the cone version but we do need the estimates on the energy and on the distortion).

\be{proof}
The argument is essentially the same as in the proof of Lemma 5.1.\ in \cite{Markovic14}. 
 \e{proof}

\subsection{Partial conformal naturality and families of good extensions}
We can now define the family of good extensions that we use to get the initial map in Theorem \ref{thm:heat}. Recall that $\mathbf{QC}_a(\set{S}^{n-1})$ denotes the set of quasiconformal maps $\set{S}^{n-1}\rightarrow \set{S}^{n-1}$ that fix the point $a\in\set{S}^{n-1}$.

\be{defn}
\label{defn:gadefn}
Let $a\in\set{S}^{n-1}$ and $f\in \mathbf{QC}_a(\set{S}^{n-1})$.
We identify $\set{H}^n\equiv \mathbf{H}^n$ and $\set{S}^{n-1}\equiv\ol{\set{R}^{n-1}}$ such that $a\equiv \it$. The extension $\curly{G}_a(f):\set{H}^n\rightarrow \set{H}^n$ is then defined as $\G_\it(f)$ with $\G_\it$ given by \eqref{eq:gtildedefn}.
\e{defn}
There is more than one way of identifying    $\set{H}^n\equiv \mathbf{H}^n$ and $\set{S}^{n-1}\equiv\ol{\set{R}^{n-1}}$ such that $a\equiv \it$.
An obvious question is whether the definition of $\curly{G}_a(f)$ depends on the choice of identification. But we have seen in Proposition \ref{prop:GE} (i) that  $\G_\it$ is partially conformally natural under isometries fixing infinity, and this yields that $\curly{G}_a(f)$ is well defined.

The following notion of partial conformal naturality generalises the classical notion of conformal naturality. 
\be{defn}
\label{defn:pcn}
Let $\{\curly{E}_a\}_{a\in \set{S}^{n-1}}$ be a family of extensions $\mathbf{QC}_a(\set{S}^{n-1})\rightarrow C^2(\set{H}^n)$. The family satisfies the \textbf{partial conformal naturality}, if for any two points 
$a,b\in \set{S}^{n-1}$ and any two isometries $I,J\in \mathrm{Isom}(\set{H}^n)$ with $I(b)=J(b)=a$, 
\beq
\label{eq:pcnaturality}
I \circ\curly{E}_b(f)\circ J^{-1}= \curly{E}_a(I \circ f\circ J^{-1})
\eeq
holds for all $f\in\mathbf{QC}_b(\set{S}^{n-1})$.
\e{defn}

\be{prop}
\label{prop:pcn}
The family $\{\curly{G}_a\}_{a\in \set{S}^{n-1}}$ from Definition \ref{defn:gadefn} satisfies partial conformal naturality.
\e{prop}

\be{proof}
This follows directly from the partial conformal naturality of $\G_\it$ under $\mathrm{Isom}_\it(\mathbf{H}^n)$ that was noted in Proposition \ref{prop:GE} (i).
\e{proof}


%

\subsection{Admissibility}

Next, we formulate what it means for a family of extensions (indexed by  boundary points) to be  \emph{admissible}, compare Definition 3.1 in \cite{Markovic14}.

\be{defn}[Admissible family]
\label{defn:admissible}
We say a family of extensions $\{\curly{E}_a\}_{a\in \set{S}^{n-1}}$ with $\curly{E}_a:\mathbf{QC}_a(\set{S}^{n-1})\rightarrow C^2(\set{H}^n)$ is \emph{admissible} if it satisfies the following properties.
\be{enumerate}[label=(\roman*)]
\item \textbf{Uniform quasi-isometry:} There exist constants $L=L(K)$ and $A=A(K)$ such that for every $a\in\set{S}^{n-1}$ and every $K$-qc $f\in \mathbf{QC}_a(\set{S}^{n-1})$, $\curly{E}_a(f)$ is an $(L,A)$-quasi-isometry.
\item \textbf{Uniformly bounded tension:} There exists a constant $T=T(K)>0$ such that for every $a\in\set{S}^{n-1}$ and every $K$-qc $f\in \mathbf{QC}_a(\set{S}^{n-1})$,
\beqs
\|\tau(\curly{E}_a(f))\|\leq T
\eeqs
\item \textbf{Continuity in $f$ and $a$:} Assume the sequence of $K$-qc maps $f_k\in \mathbf{QC}_{ a_k}(\set{S}^{n-1})$ converges pointwise to some $K$-qc map $f:\set{S}^{n-1}\rightarrow \set{S}^{n-1}$ and $ a_k\rightarrow  a$. Then, $f(a)=a$ and $\G_{ a_k}(f_k)\rightarrow \G_{ a}(f)$ in $C^2$-sense (i.e.\ first and second derivatives converge to those of $\G_{ a}(f)$, uniformly on compacts).
\e{enumerate}
\e{defn}

We have

\be{prop}
\label{prop:admissible}
The family $\{\curly{G}_a\}_{a\in \set{S}^{n-1}}$ from Definition \ref{defn:gadefn} is admissible in the sense of Definition \ref{defn:admissible}.
\e{prop}

\be{rmk}
Definition \ref{defn:admissible} is the analogue of Definition 3.1 in \cite{Markovic14} of an admissible extension. Notice however that the continuity of the entire family in $f$ and $a$ as stated in (iii) above is a stronger statement than the continuity of each individual $\curly{G}_a$ in $f$. We will use this stronger version in the proof of the Sector Lemma \ref{lm:sector}.
\e{rmk}

\be{proof}
The proofs of (i) and (ii) are essentially the same as for $\G_\it$ \cite{Markovic14}. We emphasize that the constants $L,A,T$ do not depend on $a\in\set{S}^{n-1}$ because all the $\G_a$ satisfy partial conformal naturality (in particular they are related by isometries). Moreover, note that we can normalize any sequence of $K$-qc maps (to get compactness) by composing with appropriate isometries and again using partial conformal naturality.

We come to the proof of (iii). The first part, $f(a)=a$, follows easily from the uniform H\"older continuity of $K$-qc maps. Indeed, $|a-f(a)|\leq |a-a_k|+|f_k(a_k)-f_k(a)|+|f_k(a)-f(a)|\rightarrow 0$, where the middle term vanishes by uniform H\"older continuity and the convergence $a_k\rightarrow a$. For the second part, take $I_k\in\mathrm{Isom}(\set{H}^n)$ with $I_k(a_k)=a$ and such that $I_k\rightarrow \mathrm{Id}$ in $C^2$ sense, uniformly on compacts, as $k\rightarrow\it$ (such $I_k$ because $a_k\rightarrow a$). By partial conformal naturality, we have
\beq
\label{eq:tcgdifference}
\G_{a_k}(f_k)-\G_a(f)= J^{-1}_k \ci \G_{a}(J_k\ci f_k\ci J_k^{-1})\ci J_k - \G_a(f)
\eeq
Since each $\G_a$ is continuous in $f$ uniformly on compacts (see Definition 3.1 in \cite{Markovic14} and recall that $\G_a$ is related to $\G_\it$ via isometries), we conclude that
\beqs
\G(J_k\ci f_k\ci J_k^{-1})\rightarrow \G_a(f),
\eeqs
holds in $C^2$-sense, uniformly on compacts, as $k\rightarrow\it$. This convergence is preserved under composition and so we find that \eqref{eq:tcgdifference} and its derivatives converge to zero, uniformly on compacts. This finishes the proof of admissibility.
\e{proof}


\section{Proof of the Sector Lemma}
\label{sect:sp}
The proof will be by contradiction. Assuming that there exists a ``bad'' admissible sector for large enough $\rho_{\mathrm{min}}$, one can bring it into a nice shape by using appropriate isometries (this is possible because of the scaling factor $e^{-\rho_{\mathrm{min}}}$ in Definition \ref{defn:sadmissible} of admissible sectors). From compactness of the set of uniform quasi-isometries fixing a point (a version of the Arzela-Ascoli theorem), partial conformal naturality of the good extension (see Definition \ref{defn:pcn}) and the fact that the tension field of the good extension is small at a ``random'' point (see Proposition \ref{prop:cone}), one then gets a contradiction.


\subsection{The contradiction assumption}
Suppose the claim is false. That is, suppose there exist $\al_0>1$, $\de_0>0$ and sequences of points $a_k\in\set{S}^{n-1}$, of $K$-qc maps $f_k\in\mathbf{QC}_{a_k}(\set{S}^{n-1})$, of numbers $\rho_k\geq k$ of $(\al_0,\rho_k)$-admissible sets $\Om_k\subset\set{S}^{n-1}$ (in the sense of Definition \ref{defn:sadmissible}) and of points $x_k\in \set{H}^n$ such that
\beq
\label{eq:forcontra}
\frac{\int\limits_{\rho_k}^{\rho_k+r_1} \int\limits_{\Om_k} |\tau(\curly{G}_{a_k}(f_k))(\rho,\zeta)|^2 \d \zeta\, \d \rho}{\int\limits_{\rho_k}^{\rho_k+r_1} \int\limits_{\Om_k} 1\, \d \zeta\, \d \rho} \geq \de_0.
\eeq
holds for all $1\leq r_1\leq k$. Here $(\rho,\zeta)$ denotes the geodesic polar coordinates centered at $x_k$. We will eventually get a contradiction to \eqref{eq:forcontra} by proving that the left hand side can be made arbitrarily small as $k\rightarrow\it$ and $r_1\rightarrow\it$. To take the limit in $k$ we need two things: convergence of the tension field (via compactness) and convergence of the geodesic polar coordinates (to the horocyclic  coordinates).

\subsection{The upper half space model of hyperbolic space}
\begin{figure}[t]
\be{center}
\includegraphics[height=8cm]{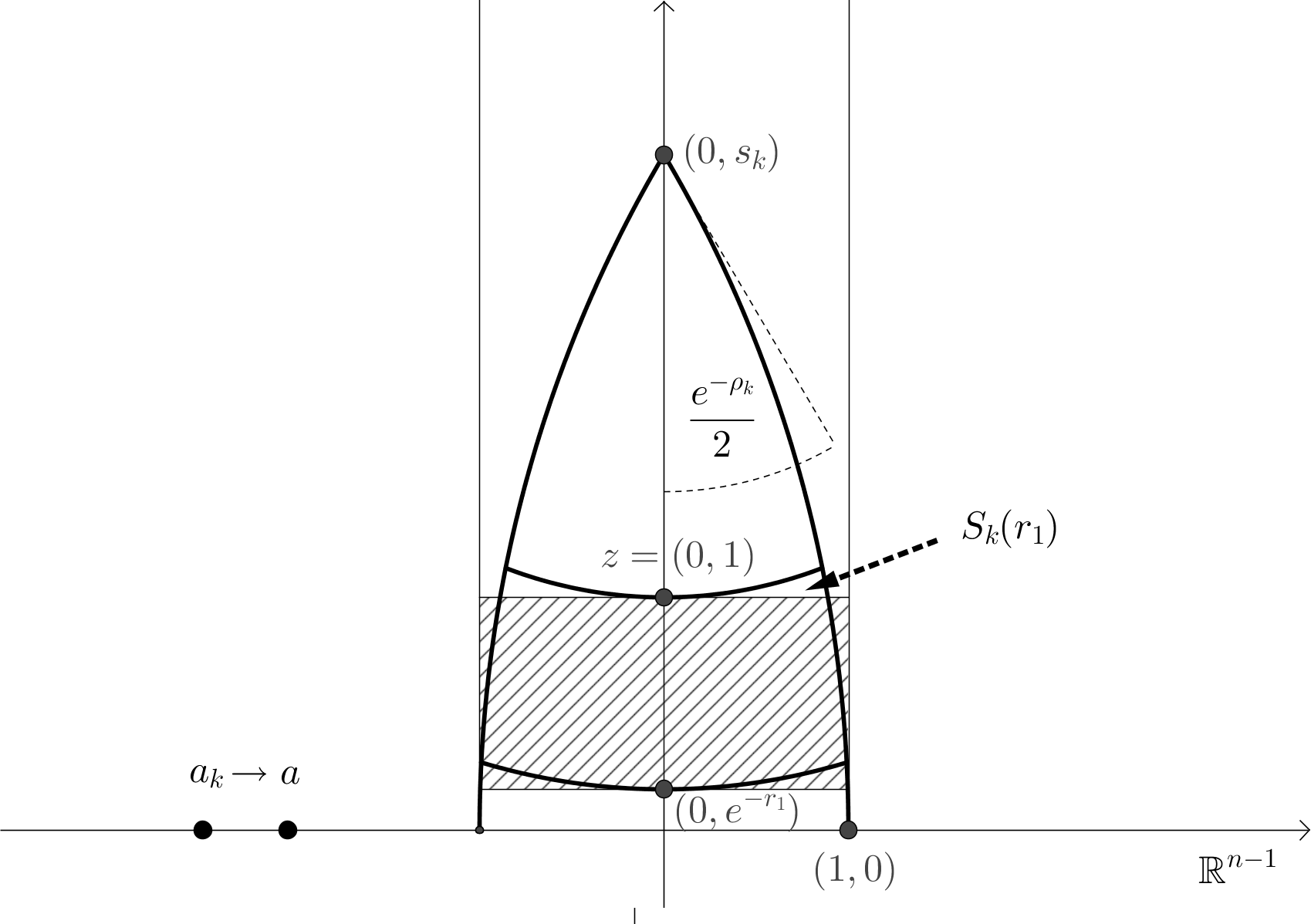}
\caption{
This picture shows how the geodesic polar coordinates centered at $x_k=(0,s_k)$ converge to horocyclic coordinates as $k\rightarrow\it$.
We see a cross cut of the upper half space model, all Euclidean coordinates $(b,s)$ should be read as $(0,\ldots,0,b,s)\in\mathbf{H}^n$.
It is intuitively clear that, as $s_k\rightarrow\it$, the cross cut of the sector $S_k(r_1)$ will ``flatten out'' and converge to the shaded region (our proof only uses the containments expressed as \eqref{eq:contain}). Notice that the geodesic which makes an initial ``angle'' with the $(0,1)$ axis of order $e^{-\rho_k}$ ends at a boundary point which is of order one as $k\rightarrow\it$.}
\label{fig:UHP}
\e{center}
\end{figure}

We work in the upper half space model of hyperbolic space $\mathbf{H}^n$. We call $D_{\n,k}$, $D_{\out,k}$ the disks which exist by Definition \ref{defn:sadmissible} since $\Om_k$ is $(\al_0,\rho_k)$-admissible. Without loss of generality, we may assume that the disks have the same center, call it $c_k\in\set{S}^{n-1}$ (otherwise this can be achieved by changing $\al_0$ to $2\al_0$). We identify $\set{H}^n\equiv \mathbf{H}^n$  such that
\beqs
x_k\equiv (0,\ldots,0,s_k),\quad (\rho_k, c_k) \equiv z=(0,\ldots,0,1),
\eeqs
(such an identification is not unique). See Figure \ref{fig:UHP} for a picture of the situation. Here $s_k>0$ is determined by the condition $$d_{\set{H}^n}((0,\ldots,0,s_k),(0,\ldots,0,1))=\rho_k.$$  It is helpful in the following to keep in mind that $s_k\sim 4e^{\rho_k}\rightarrow\it$ as $\rho_k\rightarrow\it$ (the notation $\sim$ means that $\lim_{k\rightarrow\it}\frac{s_k}{4e^{\rho_k}}= 1$).

We identify $a_k,f_k$ with their realizations in the upper half space model, $a_k\in\ol{\set{R}^{n-1}}$ and $f_k\in \mathbf{QC}_{a_k}(\ol{\set{R}^{n-1}})$. (The reader may be surprised that we do not require $a_k=\it$ in the upper half space model and instead choose an upper half space model that gives $x_k,(\rho_k, c_k)$ the nice coordinates above. The reason is that this chart is well fitting to see the convergence of the geodesic polar coordinates to horocyclic coordinates as $\rho_k\rightarrow\it$.)

We post-compose $f_k$ by a sequence of isometries such that the resulting sequence fixes a point \emph{inside} $\mathbf{H}^n$. That is, we find $I_k\in \mathrm{Isom}_{a_k}(\mathbf{H}^n)$ such that
\beqs
I_k (\curly{G}_{a_k}(f_k)(z))=z
\eeqs
and define
\beqs
g_k:=I_k\ci f_k.
\eeqs
which then satisfies $g_k(z)=z$. Note also that $g_k\in\mathbf{QC}_{a_k}(\ol{\set{R}^{n-1}})$ and so by the partial conformal naturality of the good extension (in the sense of Definition \ref{defn:pcn})
\beqs
\G_{a_k}(g_k)=I_k \ci \G_{a_k}(f_k),
\eeqs
and in particular
\beq
\label{eq:Geq}
|\tau(\G_{a_k}(f_k))|=|\tau(\G_{a_k}(g_k))|.
\eeq 

\subsection{Convergence of the tension from compactness}
Since $\ol{\set{R}^{n-1}}$ is compact,  up to passing to a subsequence, there exists $a\in \ol{\set{R}^{n-1}}$ such that $a_k\rightarrow a$. 

Moreover, by using standard arguments about quasiconformal maps and quasi-isometries (in particular an extension of the Arzela-Ascoli theorem for uniform quasi-isometries which all fix the same point) one proves that, up to passing to a subsequence, there exists a quasiconformal map $g: \ol{\set{R}^{n-1}} \to \ol{\set{R}^{n-1}}$ such that $g_k\rightarrow g$ pointwise. 

Together, these facts enable us to apply the continuity of the good extension in the sense of Definition \ref{defn:admissible} (iii). First, this implies $g(a)=a$ and so $g\in\mathbf{QC}_{a}(\ol{\set{R}^{n-1}})$. Second, it implies that $\G_{a_k}(g_k)\rightarrow \G_{a}(g)$ in $C^2$-sense, uniformly on compacts. The upshot of this first part of the proof is that we have
\beq
\label{eq:tauconv}
|\tau(\curly{G}_{a_k}(g_k))|\rightarrow |\tau(\G_{a}(g))|
\eeq
pointwise, uniformly on compacts.

\subsection{Convergence of geodesic polar coordinates to horocyclic coordinates}
For this part, it is helpful to consider Figure \ref{fig:UHP}.  Let $(b,h)\in \set{R}^{n-1} \times \set{R}$ denote horocyclic coordinates on the upper half-space $\mathbf{H}^n$.  Recall that a point in $\mathbf{H}^n$ with coordinates $(b,h)$ lies above $b \in \set{R}^{n-1}$ and the horosphere through this point has the signed hyperbolic  height $h$ (with the normalization that the horosphere through the point $z=(0,\ldots,1)$ has the height $0$). Note that the point in $\mathbf{H}^n$ with  horocyclic coordinates $(b,h)$ has the Euclidean coordinates $(b,e^{-h})$.

Given a point in geodesic polar coordinates $(\rho,\zeta)$, we identify $\zeta\in\set{S}^{n-1}$ with the ``endpoint'' of the corresponding geodesic in $\ol{\set{R}^{n-1}}$. (More precisely, we recall that $\zeta$ is identified with an element of the unit tangent space at the point where the geodesic polar coordinates are centered, here $x_k$. Then, we find the endpoint of the geodesic with this unit tangent vector as the initial direction and call this endpoint $\zeta$.) With this identification, we have $\zeta \to b$ and $(\rho-\rho_k) \to h$, when $k\to \infty$.

Recall that the integration in \eqref{eq:forcontra} takes place over the sector
\beqs
S_k(r_1)\equiv S_k(x_k,\rho_k,r_1,\Om_k)=\setof{(\rho,\zeta)}{\rho_k\leq \rho \leq \rho_k+r_1,\, \zeta\in \Om_k}.
\eeqs
Using \eqref{eq:Geq}, we may rewrite \eqref{eq:forcontra} as
\beq
\label{eq:integrand}
\de_0 \leq \frac{\int\limits_{S_k(r_1)} |\tau(\curly{G}_{a_k}(g_k))(\rho,\zeta)|^2 \d\rho\, \d\zeta}{\int\limits_{S_k(r_1)} 1\, \d\rho\, \d\zeta},
\eeq
where $\d\zeta$ is the measure on $\ol{\set{R}^{n-1}}$ induced by the spherical measure. We now discuss the limiting properties, as $k\rightarrow\it$, of $S_k(r_1)$ where $\Om_k$ is identified with an appropriate subset of $\ol{\set{R}^{n-1}}$ in the way discussed above. To this end, we define the cylinders
\beq
\label{eq:cyldefn}
\begin{aligned}
\mathrm{Cyl}_{\mathrm{in}}(r_1):&=\setof{(b,h)\in\set{R}^{n-1}\times \set{R} }{|b|\leq 2 \al_0^{-1},\, 0\leq h\leq r_1}\\
\mathrm{Cyl}_{\mathrm{out}}(r_1):&=\setof{(b,h)\in\set{R}^{n-1}\times \set{R} }{|b|\leq 2\al_0,\, 0\leq h\leq r_1}.
\end{aligned}
\eeq
We recall that $\Om_k$ viewed as a subset of $\set{S}^{n-1}$ is $(\al_0,\rho_k)$-admissible in the sense of Definition \ref{defn:sadmissible}. That is, there exist disks $D_{\n,k},D_{\out,k}\subset\set{S}^{n-1}$ such that
\beq
\label{eq:admissiblerec}
D_{\n,k}\subset\Om_k\subset D_{\out,k}
\eeq
and the radius of $D_{\n,k}$ is at least $\al_0^{-1} e^{-\rho_k}$, while the radius of $D_{\out,k}$ is at most $\al_0e^{-\rho_k}$. We also recall that we assumed that $D_{\n,k},D_{\out,k}$ are centered at the same point $c_k\in\set{S}^{n-1}$ which is identified with the downward pointing normal in our upper half space model of hyperbolic space. 

When we identify $\Om_k$ with a subset of $\ol{\set{R}^{n-1}}$ as discussed above, \eqref{eq:admissiblerec} yields
\beq
\label{eq:contain}
\limsup_{k\rightarrow \it}\ind_{S_k(r_1)}\leq \ind_{\mathrm{Cyl}_{\mathrm{out}}(r_1)},\quad 
\liminf_{k\rightarrow \it}\ind_{S_k(r_1)}\geq \ind_{\mathrm{Cyl}_{\mathrm{in}}(r_1)}
\eeq
for all $1\leq r<\it$. Here $\ind$ denotes the characteristic function of a subset of hyperbolic space. The reader may find it helpful to consider Figure \ref{fig:UHP}. (The relations \eqref{eq:contain} together with the fact that the bounds on the $b$ variable in \eqref{eq:cyldefn} are independent of $k$ are the manifestations of admissible sectors having ``bounded geometry'' near the boundary. Note that the factor $e^{-\rho_k}$ in Definition \ref{defn:sadmissible} of an $(\al_0,\rho_k)$-admissible $\Om_k$ is important for this.) 

For a fixed $r_1\geq 1$, we take $k\rightarrow \it$ in \eqref{eq:integrand}, more precisely we take the $\limsup_{k\rightarrow \it}$ of the numerator and the $\liminf_{k\rightarrow \it}$ of the denominator in \eqref{eq:integrand}. It is elementary to check that $\d\rho \to \d h$, when $k \to \infty$, and 
\beq
\lim_{\rho_k\rightarrow\it} 2e^{\rho_k} \d\zeta=\d b,
\eeq
where $\d b$ is the standard Lebesgue measure on $\set{R}^{n-1}$. Recall also \eqref{eq:tauconv} which says that $|\tau(\curly{G}_{a_k}(g_k))|\rightarrow |\tau(\G_{a}(g))|$ pointwise as $k\rightarrow\it$. We can then use dominated convergence together with the relations \eqref{eq:contain} to conclude from \eqref{eq:integrand} that
\beq
\label{eq:cylindrical}
\begin{aligned}
\de_0 &\leq \frac{\int\limits_{\mathrm{Cyl}_{\mathrm{out}}(r_1)} |\tau(\curly{G}_{a}(g))(b,h)|^2 \d b\, \d\zeta}{\int\limits_{\mathrm{Cyl}_{\mathrm{in}}(r_1)} 1\, \d b\, \d h}\\
&=\frac{C(\al_0)}{r_1} \int\limits_{\mathrm{Cyl}_{\mathrm{out}}(r_1)} |\tau(\curly{G}_{a}(g))(b,h)|^2 \d b\, \d h\\
&=\frac{C(\al_0)}{r_1} \int\limits_0^{r_1} \int\limits_{|b|\leq 2\al_0}|\tau(\curly{G}_{a}(g))(b,h)|^2 \d b\, \d h
\end{aligned}
\eeq
holds for all $1\leq r_1<\it$. Here $C(\al_0)>0$ is an appropriate constant.

\subsection{Getting a contradiction}
\be{lm}
\label{lm:forcontra}
For almost every $b\in\set{R}^{n-1}$, we have
\beq
\label{eq:pwiseconv}
\lim_{r_1\rightarrow\it}\frac{1}{r_1} \int\limits_0^{r_1} |\tau(\curly{G}_{a}(g))(b,h)|^2  \d h=0.
\eeq
\e{lm}

By dominated convergence, Lemma \ref{lm:forcontra} gives a contradiction to \eqref{eq:cylindrical}. To prove the Sector Lemma, it therefore remains to give the

\be{proof}[Proof of Lemma \ref{lm:forcontra}]
 Let $\de'>0$. The lemma will follow easily once we prove the following \emph{claim}: For almost every $b\in\set{R}^{n-1}$, there exists $r_2=r_2(f,b,\de')$ such that for all $s\geq r_2$,
\beq
\label{eq:de'}
|\tau(\G_a(g))(b,h)|^2<\de'.
\eeq
By Proposition \ref{prop:cone} we know that  $|\tau(\G_a(g))(w)|^2<\de'$, when $w \to b$ and $w$ belongs to the cone around the geodesic connecting $a$ and $b$ (the cone contains all the points that are within some fixed distance from the geodesic connecting $a$ and $b$).  But, any geodesic converging to $b$ will eventually enter this cone, and so will the geodesic starting at $\infty$ and ending at $b$. This proves the claim.

We let $r_1>r_2$. We can now cut the integral from \eqref{eq:pwiseconv} into a bad part (where we use that $\|\tau(\G_a(g))\|\leq T$) and a good part (where \eqref{eq:de'} holds) :
\beqs
\begin{aligned}
&\frac{1}{r_1} \int\limits_0^{r_1} |\tau(\G_a(g))(b,h)|^2 \d h\\
&= \frac{1}{r_1} \int\limits_0^{r_2} |\tau(\G_a(g))(b,h)|^2 \d h + \frac{1}{r_1} \int\limits_{r_2}^{r_1} |\tau(\G_a(g))(b,h)|^2 \d h\\
&\leq T^2 \frac{r_2}{r_1} +\de'.
\end{aligned}
\eeqs
The first term vanishes as $r_1\rightarrow\it$. Since $\de'>0$ was arbitrary, this proves \eqref{eq:pwiseconv}.
\e{proof}

\section{Proof of Theorem \ref{thm:close}}
As mentioned before, the proof is a straightforward  generalization of the arguments in \cite{Markovic14} to higher dimensions and the observation that the estimates have enough ``wiggle room'' to allow for a sufficiently small $\eps_0$. Consequently, we only give a sketch of the argument here and refer the reader to \cite{Markovic14} for a more thorough discussion. 

We work in the unit ball model of hyperbolic space which we denote by $\set{B}^n$.  
Let $f\in\mathbf{QC}_a(\set{S}^{n-1}$) be a $K$-qc map and let $\psi:\set{B}^n\rightarrow\set{B}^n$ be a $C^2$ quasi-isometry with the boundary map $f$. Then
\beqs
\|d_{\set{H}^{n}}\l(\curly{G}_a(f),\psi\r)\|<\it
\eeqs
since both maps are quasi-isometries which extend $f$. As in \cite{Markovic14}, we may assume without loss of generality \ that
\beq
\label{eq:d0}
\bd(f)(0)\geq \|\bd(f)\|-D-1
\eeq
where $D=D(K)$ and 
\beqs
\bd(f)(x)\equiv d_{\set{H}^n}\l(\curly{G}_a(f)(x),\psi(x)\r).
\eeqs


As in  \cite{Markovic14}, combining  \eqref{eq:d0} and Green's identity for $\bd^2(f)$ we obtain the crucial estimate
\beq
\label{eq:LHS}
\int\limits_{\set{B}^n} \bg_r(x) \Delta \bd^2(f)(x)\d\blam\leq D'\|\bd(f)\| +D''.
\eeq
Here $D'=2(D-1)$ and $D''=(D+1)^2$ depend only on $K$ and $\d\blam$ is the hyperbolic volume measure in the unit ball model. In this section only, $x=(\rho,\zeta)$ stands for Euclidean polar coordinates, i.e.\ $\rho\in[0,1)$, and \emph{not} for the geodesic polar coordinates (we do this for the sake of comparability with \cite{Markovic14}). 
\beq
\label{eq:volumemeasure}
\d\blam(x)=\frac{n\rho^{n-1}}{(1-\rho^2)^n} \d\rho\,\d\sigma(\zeta)
\eeq
where $\d\sigma$ is the Lebesgue measure on $\set{S}^{n-1}$, normalized to $\sigma(\set{S}^{n-1})=1$. Finally, $\bg_r$ is the Green's function of $-\Delta$ on $r\set{B}^n,\,0\leq r<1$. Explicitly \cite{Ahlfors},
\beqs
\begin{aligned}
\bg_r(x)=
\frac{1}{n}\int\limits_{|x|}^r \frac{(1-s^2)^{n-2}}{s^{n-1}} \d s,\quad |x|\leq r
\end{aligned}
\eeqs
and $\bg_r(x)=0,\,r<|x|<1$. Note that $\bg_r$ is a radial function. We often abuse notation and write $\bg_r(\rho)$ for $\rho>0$. We have the lower bound 
\beq
\label{eq:GFestimate}
\bg_r(\rho)\geq C_\bg \frac{(1-\rho^2)^{n-1}}{\rho^{n-2}},
\eeq
where $C_\bg$ is a universal constant. Moreover, $\bg_r\rightarrow \bg_1$ uniformly on compacts as $r\rightarrow 1$.

Consider \eqref{eq:LHS}. Note that the claim that $\|\bd\|$ is bounded by a constant would follow if we had a lower bound on the left hand side of the form $(D'+1)\|\bd(f)\|$. This is what is done in \cite{Markovic14}, and the same proof can be repeated word by word modulo two minor modifications (one in Lemma 3.2 and one in Lemma 4.2 from 
\cite{Markovic14}) which we describe below.

The next step in \cite{Markovic14} is to estimate the set where $\bd^2(f)$ is small, see Lemma 4.1,  and the proof generalizes directly to higher dimensions.  
An important tool in the proof of the main Lemma 4.2 in \cite{Markovic14} are the following estimates from  \cite{TamWan98}, originally from \cite{SchoenYau79, JagerKaul}. They say that for any $F,G\in C^2(\set{H}^n)$
\beq
\label{eq:tamwan1}
\Delta\bd^2\geq -2\bd \l(\|\tau(F)\|+\|\tau(G)\|\r),\quad \text{on } \set{H}^n,
\eeq
where $\bd\equiv d_{\set{H}^n}(F,G)$. Moreover, for all $K_1\geq 1$ there exists $q=q(K_1)>0$ such that
\beq
\label{eq:tamwan2}
\Delta\bd^2\geq -2\bd \l(\tau(F)+\tau(G)\r)+2q\, \bd\, \mathbf{e}(F) \tanh\l(\frac{\bd}{2}\r),
\eeq
holds for all $x\in\set{H}^n$ with $\mathbf{K}(F)(x)\leq K_1$. One follows the proof of Lemma 4.2 in \cite{Markovic14} and applies these estimates. The only difference is that one takes
\beq
\label{eq:eps0defn}
\eps_0(K):=\frac{q(2K)}{8}\tanh(1/4),
\eeq
which has a relative factor of $1/2$ compared to the definition of $\eps_0(K)$ on page 19 of \cite{Markovic14}. This is exactly the place where we use that we have some ``wiggle room''.

The second place where justification is required is to show that Lemma 3.2 from \cite{Markovic14} holds in $n$ dimensions. The argument from \cite{Markovic14} applies provided that
\beq
\label{eq:ffixed}
\lim_{r\rightarrow 1} \int\limits_{\set{B}^n} \bg_r(x)  \d \blam(x)=\it.
\eeq
To see this, we express the previous integral in the Euclidean  polar coordinates and use \eqref{eq:volumemeasure} and \eqref{eq:GFestimate} to get
\beqs
\int\limits_{0}^{r} \frac{n\rho^{n-1}\bg_r(\rho)}{(1-\rho^2)^{n-1}} \d \rho  \geq \frac{n C_\bg}{2} \log\l(\frac{1}{1-r^2}\r)\rightarrow \it,\quad r\rightarrow 1.
\eeqs
The analogue of Lemma 3.2 then follows by the usual compactness argument (we can pre- and postcompose by appropriate isometries to normalize the $K$-qc maps thanks to the partial conformal naturality of the good extension). 

\be{appendix}

\section{Heat travels ballistically in hyperbolic space}
\label{app:heat}
In this appendix, we discuss the diffusion of heat in hyperbolic space. It is known that heat travels approximately ballistically in the hyperbolic space. By this we mean that, for large $t$, the measure whose density is given by the heat kernel $H(x,y,t)$ times the hyperbolic volume measure $\d\lam(y)$ is effectively supported on a certain ``main annulus''(in geodesic polar coordinates), which is centered at $x$ and has inner and outer radii of order $t$ (see e.g.\ Corollary 5.7.3.\ in \cite{Davies}). (In Euclidean space, such an annulus would have radii of order $\sqrt{t}$.)

Here we prove a more precise version. It says that the main annulus has $\rho$-values of the form $(n-1)t\pm r\sqrt{t}$ with $r=O(1)$ distributed according to the standard Gaussian measure $e^{-r^2/4}\d r$ on the main annulus, see Figure \ref{fig:Gaussian} for a picture.

While these facts are presumably known to experts, we could not find a reference. Therefore we discuss this topic here in some detail. The proof only uses the heat kernel bounds in \cite{DaviesMandouvalos87}.

\be{figure}[t]
\be{center}
\includegraphics[height=7.5cm]{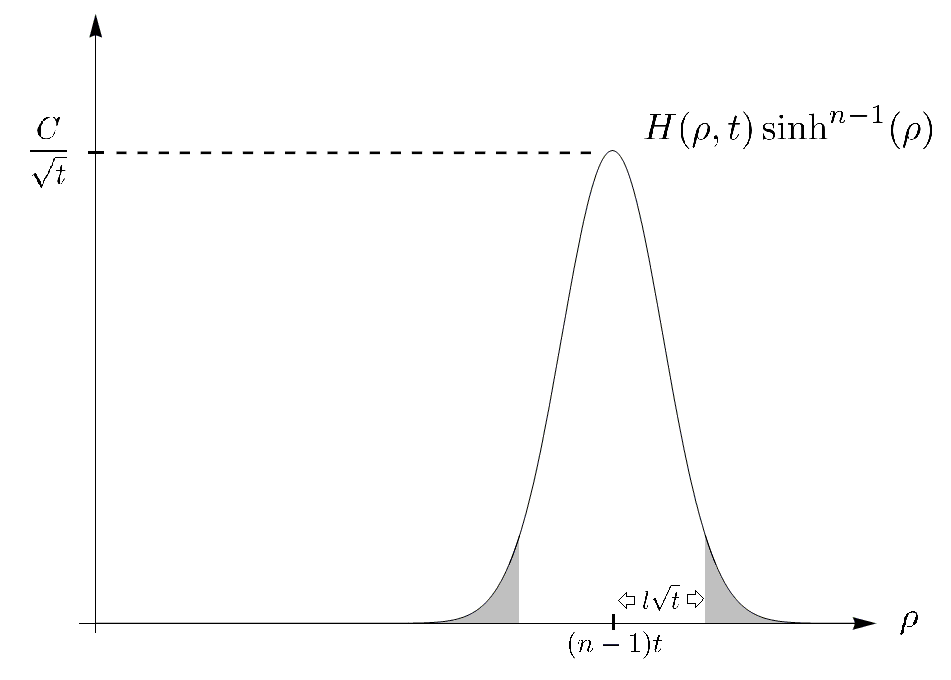}
\caption{This plot of the heat kernel times hyperbolic volume measure as a function of the radial coordinate $\rho$ shows how heat is transported in hyperbolic space (for large $t$). The function is centered at $(n-1)t$ and decays around that center on scale $\sqrt{t}$ like the standard Gaussian. The choice $l=l(\eps)$ from Proposition \ref{prop:annulus} (i) is such that the shaded region has area $\eps$. Thus, the function $H(\rho,t)\sinh^{n-1}(\rho)$ is mainly supported on the region in between and this defines the ``main annulus''.}
\label{fig:Gaussian}
\e{center}
\e{figure}

\be{prop}
\label{prop:annulus}
\be{enumerate}[label=(\roman*)]
\item
There exists a universal (dimension dependent) constant $C_n>0$, such that for all $\eps>0$ and all $t\geq 1$,
\beq
 \int\limits_{|\rho-(n-1)t|> l(\eps)\sqrt{t}} H(\rho,t)\sinh^{n-1}(\rho) \d \rho < \eps.
\eeq
where $l(\eps):=\sqrt{8\log\l(\frac{C_n}{\eps}\r)}$.

\item Let $\Phi:\set{R}_+\rightarrow \set{R}_+$ be a bounded measurable function. Let $l\geq 1$. Then, for all $t\geq 2l^2$,
\beq
\label{eq:annulusGaussian}
\begin{aligned}
&\int\limits_{|\rho-(n-1)t|\leq l\sqrt{t}} \Phi(\rho) H(\rho,t)\sinh^{n-1}(\rho)  \d \rho\\
&\be{cases}
\geq \frac{1}{C_n'}\int\limits_{-l}^l  \Phi((n-1)t+r\sqrt{t})\, e^{-\frac{r^2}{4}}  \d r\\
\leq C_n'\int\limits_{-l}^l  \Phi((n-1)t+r\sqrt{t})\, e^{-\frac{r^2}{4}}  \d r.
\e{cases}
\end{aligned}
\eeq
where $C_n'>1$ is a universal (dimension dependent) constant.

\item Let $\Phi:\set{R}_+\rightarrow \set{R}_+$ be a bounded measurable function. Let $l(\eps)$ be as in (i) and let $C_n'$ be the universal constant from (ii). Then, for all $\eps>0$ and all $t\geq 2l(\eps)^2$,
\beq
\label{eq:see2}
\begin{aligned}
\int\limits_{0}^\it \Phi(\rho) H(\rho,t)\sinh^{n-1}(\rho)\d \rho \leq \frac{C_n'}{\sqrt{t}} \int\limits_{(n-1)t-l(\eps)\sqrt{t}}^{(n-1)t+l(\eps)\sqrt{t}} \Phi(\rho) \d \rho + \eps
\end{aligned}
\eeq
\e{enumerate}
\e{prop}

\be{rmk}
This proposition is used to prove Theorem \ref{thm:heat}. That proof would hold under weaker assumptions on the form of heat diffusion (for instance, it would be enough to know that the effective support of $H(x,y,t)\d\lam(y)$ is an annulus centered at $x$ which has inner radius going to infinity and diverging width as $t\rightarrow\infty$). Nonetheless, we give a precise description of the heat diffusion because this may be of independent interest and because the proof is straightforward.
\e{rmk}

Statement (iii) follows directly from (i) and the upper bound in (ii). In the main text, we apply statement (iii) to $\Phi$ being the spherical average of $|\tau(\curly{G}_a(f))|^2$, see \eqref{eq:see}. We will not use the lower bound in (ii), it is only stated here for the sake of completeness. 

\be{proof}
Throughout the proof, we write $C>0$ for a universal (dimension dependent) constant; the numerical value of $C$ may change even in the same line. We first prove statement (i). By Theorem 3.1 in \cite{DaviesMandouvalos87},
\beqs
H(\rho,t)\leq  C t^{-n/2} (1+\rho+t)^{\frac{n-3}{2}}(1+\rho)  \exp\l(-\frac{\rho^2}{4t}-\frac{(n-1)^2}{4}t-\frac{n-1}{2}\rho\r).
\eeqs
Since $\sinh(\rho)< \frac{\exp(\rho)}{2}$ for $\rho>0$, we get
\beq
\begin{aligned}
\label{eq:recall}
&\sinh^{n-1}(\rho) H(\rho,t) <\\
& C t^{-n/2} (1+\rho+t)^{\frac{n-3}{2}} (1+\rho) \exp\l(-\frac{1}{4}\l(\frac{\rho}{\sqrt{t}}-(n-1)\sqrt{t}\r)^2\r).
\end{aligned}
\eeq
We change variables to $r=\frac{\rho}{\sqrt{t}}-(n-1)\sqrt{t}$ and find, for all $t\geq 1$ and $l>0$ to be determined,
\beqs
\begin{aligned}
&\int\limits_{|\rho-(n-1)t|> l\sqrt{t}}\sinh^{n-1}(\rho) H(\rho,t)\d \rho \leq\\
&  C \int\limits_{|r|\geq l} \l(n+1+\frac{r}{\sqrt{t}}\r)^{\frac{n-1}{2}} e^{-\frac{r^2}{4}} \d r.
\end{aligned}
\eeqs
Notice that for all $t\geq 1$
\beqs
\begin{aligned}
 C \int\limits_{|r|\geq l}\l(n+1+\frac{r}{\sqrt{t}}\r)^{\frac{n-1}{2}} e^{-\frac{r^2}{4}} \d r &\leq C e^{-\frac{l^2}{8}} \int\limits_{\set{R}} \l(n+1+r\r)^{\frac{n-1}{2}} e^{-\frac{r^2}{8}} \d r\\
 &\equiv C_n e^{-\frac{l^2}{8}}
\end{aligned}
\eeqs
where $C_n$ is defined by the last equality. Let $\eps>0$. Setting $l=l(\eps)=\sqrt{8\log\l(\frac{C_n}{\eps}\r)}$, yields
\beqs
\int\limits_{|\rho-(n-1)t|> l(\eps)\sqrt{t}}\sinh^{n-1}(\rho) H(\rho,t)\d \rho < C_n e^{-\frac{l(\eps)^2}{8}} = \eps.
\eeqs
This proves (i).

We come to statement (ii). Fix $l\geq 1$. Recall \eqref{eq:recall} and integrate it over the interior of the main annulus now. Changing variables again to $r=\frac{\rho}{\sqrt{t}}-(n-1)\sqrt{t}$ gives
\beq
\label{eq:annulusUB}
\begin{aligned}
&\int\limits_{|\rho-(n-1)t|\leq l\sqrt{t}} \Phi(\rho) H(\rho,t)\sinh^{n-1}(\rho)  \d \rho \leq\\
&   C \int\limits_{-l}^l  \Phi((n-1)t+r\sqrt{t}) \l(n+1+\frac{r}{\sqrt{t}}\r)^{\frac{n-1}{2}} e^{-\frac{r^2}{4}}  \d r.
\end{aligned}
\eeq
When $t\geq 2 l^2$, we can bound $n+1+\frac{r}{\sqrt{t}}\leq n+2$. This implies the upper bound in \eqref{eq:annulusGaussian} for an appropriate universal constant $C_n'$.

For the lower bound in \eqref{eq:annulusGaussian}, we use that Theorem 3.1 in \cite{DaviesMandouvalos87} also gives
\beqs
H(\rho,t)\geq C t^{-n/2} (1+\rho+t)^{\frac{n-3}{2}}(1+\rho)  \exp\l(-\frac{\rho^2}{4t}-\frac{(n-1)^2}{4}t-\frac{n-1}{2}\rho\r).
\eeqs
One can check that $\sinh(\rho)>\frac{1}{4} e^{\rho}$ holds for all $\rho$ with $|\rho-(n-1)t|\geq l\sqrt{t}$ and all $t\geq 2l^2,\, l\geq 1$. After integration and the change of variables $r=\frac{\rho}{\sqrt{t}}-(n-1)\sqrt{t}$, this yields the following analogue to \eqref{eq:annulusUB}
\beq
\label{eq:annulusLB}
\begin{aligned}
&\int\limits_{|\rho-(n-1)t|\leq l\sqrt{t}} \Phi(\rho) H(\rho,t)\sinh^{n-1}(\rho)  \d \rho \geq \\
&   C \int\limits_{-l}^l  \Phi((n-1)t+r\sqrt{t}) \l(n-1+\frac{r}{\sqrt{t}}\r)^{\frac{n-1}{2}} e^{-\frac{r^2}{4}}  \d r.
\end{aligned}
\eeq
Again, $n-1+\frac{r}{\sqrt{t}}$ can be bounded below by a uniform constant. This implies the lower bound in \eqref{eq:annulusGaussian} for an appropriate $C_n'$.

Finally, (iii) follows directly from (i) and (ii) by dropping the Gaussian and undoing the change of variables in \eqref{eq:annulusGaussian}.
\e{proof}
\e{appendix}

\bibliographystyle{amsplain}

\end{document}